%% file: main.tex
\documentclass[11pt, reqno, english]{amsart}  

\usepackage{preamble}

\begin{document}

\title{Trace of Homogeneous Fractional Sobolev Spaces on Strip-like Domains}
\author{Khunpob Sereesuchart}
\email{ksereesu@andrew.cmu.edu} 
\date{\today}

\maketitle

\begin{abstract}
In this paper, we discuss the trace operator for homogeneous fractional Sobolev spaces over infinite strip-like domains. We determine intrinsic seminorms on the trace space that allow for a bounded right inverse. The intrinsic seminorm includes two features previously used to describe the trace of homogeneous Sobolev spaces, a relation between the two disconnected components of the trace and the screened Sobolev seminorm. However, unlike its homogeneous Sobolev space equivalent, fractional Sobolev spaces require a far screened Sobolev seminorm that captures the non-local properties of fractional Sobolev spaces. We study some basic relationships between this new far screened Sobolev space with previously discussed screened Sobolev spaces. 
\end{abstract}

\tableofcontents

\subfile{sections/1.Introduction}

\subfile{sections/2.Preliminaries}

\subfile{sections/3.FlatOmega}

\subfile{sections/4.GeneralOmega}

\subfile{sections/5.Fourier_Transforms}

\section*{Acknowledgements} 
This paper is part of the author's honor thesis at Carnegie Mellon University This research was partially supported by the National Science Foundation under grants No DMS  2108784 and DMS  1714098. I would also like to thank my Master's advisor Giovanni Leoni for supporting me for the previous few years. 
\printbibliography
\end{document}

%% file: sections/1.Introduction.tex
\section{Introduction}
\subfile{1.Introduction/1.1.Introduction}

%% file: sections/1.Introduction/1.1.Introduction.tex

For $\Omega\subseteq \mathbb R^N$ open, $0<s<1$ and $1\leq p<\infty$, the homogeneous fractional Sobolev space, denoted as $\dot W^{s,p}(\Omega)$, is the space of all functions $u\in L^p_{\text{Loc}}(\Omega)$ such that the seminorm
$$|u|_{\dot W^{s,p}(\Omega)}:=\left(\int_\Omega \int_\Omega \frac{|u(x)-u(y)|^p}{\Vert x-y\Vert^{N+sp}}dydx\right)^{1/p}$$
is finite. The goal of this paper is to determine the trace space of $\dot W^{s,p}(\Omega)$, where $0<s<1$, $1\leq p<\infty$, and $sp>1$, for domains of the form
\begin{equation}\label{generalomega}
\Omega=\{(x',x_N)|\eta^-(x')<x_N<\eta^+(x')\}
\end{equation}
where $\eta^\pm:\mathbb R^{N-1}\rightarrow \mathbb R$ are Lipschitz with $\eta^-<\eta^+$. 
\par The trace of the non-homogeneous fractional Sobolev space, denoted $W^{s,p}(\Omega)$, is well studied. We define $f\in W^{s,p}(\Omega)$ if and only if $f\in L^p(\Omega)$ and $|f|_{W^{s,p}(\Omega)}<\infty$. Furthermore, we can equip the norm
$$\Vert f\Vert_{W^{s,p}(\Omega)}=\left(\Vert f\Vert_{L^p(\Omega)}^p+|f|_{\dot W^{s,p}(\Omega)}^p\right)^{1/p}.$$
We also define $W^{s,p}_0(\Omega)$ to be the closure of $C_c^\infty(\Omega)$ under the norm $\Vert\cdot \Vert_{ W^{s,p}(\Omega)}$. In theorem $6.76$ of \cite{fss}, the following theorem was proven.
\begin{theorem}
For $0<s<1$, $1\leq p<\infty$, $sp\leq 1$, and $\Omega\subset \R^N$ with Lipschitz continuous boundary, we have
$$W^{s,p}(\Omega)=W^{s,p}_0(\Omega)$$
\end{theorem}
Note that for $f\in C_c^\infty(\Omega)$, we have $\Tr(f)=0$. Therefore, if a continuous trace function exists over $W_0^{s,p}(\Omega)$, then the trace function is the zero function. With the previous theorem, any continuous trace function over $W^{s,p}(\Omega)$ is the zero function. Yet, there are functions with nonzero trace in $W^{s,p}(\Omega)$, which implies that we cannot define a continuous trace over $W^{s,p}(\Omega)$ for $sp\leq 1$. Therefore, we will only discuss the case when $sp>1$. 
\par For $sp>1$ and $\Omega\subset \R^N$ where $\Omega$ has a uniformly Lipschitz continuous boundary, a continuous trace is well-defined in theorem $9.49$ from \cite{fss}.
\begin{theorem}\label{generalset}
Let $\Omega\subset \R^N$, $N\geq 2$, open with $\partial \Omega$ uniformly Lipschitz continuous, $1<p<\infty$, $0<s<1$ and $sp>1$. Then there is $C>0$ such that
$$\Vert \Tr(u)\Vert_{W^{s-1/p,p}(\partial \Omega)} \leq C\Vert u\Vert_{ W^{s,p}(\Omega)}$$
for all $u\in W^{s,p}(\Omega)$. Furthermore, for every $g\in W^{s-1/p,p}(\partial\Omega)$, there exists a function $u\in W^{s,p}(\Omega)$ such that $\Tr(u)=g$ and
$$\Vert u\Vert_{W^{s,p}(\Omega)}\leq C\Vert g\Vert_{ W^{s-1/p,p}(\partial\Omega)}$$
\end{theorem}
\begin{remark}
We note that $\Omega$ as in (\ref{generalomega}) does not have a uniformly Lipschitz continuous boundary if $\inf_{x'\in \R^{N-1}}\eta^+(x')-\eta^-(x')=0$. This will be discussed later in my thesis \cite{ksereesuthesis}.   
\end{remark}
For more information on Sobolev spaces, see \cite{adams_f._adams_2003}, \cite{besov_ilin_nikolskiisergeim_taibleson_1979}, \cite{besov_ilin_nikolski_1978}, \cite{RSMUP_1957__27__284_0}, \cite{grisvard_2011}, \cite{mazya_2011}, \cite{necas_simader_necasovasarka_tronel_kufner_2013}, \cite{leoni_2017} and \cite{zhou_du_2010}. The trace of closed $k$-sets is studied in papers by Jonnson and Wallin \cite{jonsson_wallin_1978} and \cite{jonsson_wallin_1980}.\par
The reason the trace of non-homogeneous fractional Sobolev spaces are studied is because they have applications in partial differential equations. For example, we consider, for $\Omega$ bounded, the Dirichlet boundary problem
\begin{equation}\label{sampleproblem1}
\left\{
\begin{aligned}
-\Delta u=f \tab x\in \Omega\\
u = g \tab x\in \partial\Omega
\end{aligned}\right.
\end{equation}
It is known that for $g\in W^{k+1/2,2}(\partial \Omega)$ and $f\in W^{k-1,2}(\Omega)$, $k\in \N$, there is a unique solution $u\in W^{k+1,2}(\Omega)$, as discussed in \cite{taylor_2011}. Similar relations can be found in fractional Sobolev spaces. For example, we can consider the problem

\begin{equation}\label{sampleproblem2}
\left\{
\begin{aligned}
(-\Delta)^s u=f \tab x\in \Omega\\
u = g \tab x\in \R^N\setminus\Omega 
\end{aligned}\right.
\end{equation}
where $(-\Delta)^s$ denotes the fractional Laplacian of order $s$ and $\Omega$ bounded. For $f\in L^2$, $g\in \dot H^s(\R^N)$, then $u\in g+H^s_{00}(\Omega)$, as discussed in \cite{fss}. For more information on the fractional Laplacian, see \cite{MR2354493}, \cite{MR2367025} and \cite{MR2494809}. However, due to Poincare's inequality failing for $\Omega$ unbounded, solutions to (\ref{sampleproblem1}) and (\ref{sampleproblem2}) need not be in $L^p(\Omega)$, requiring us to approach the problem differently. One method is to use weighted Sobolev spaces, as discussed in \cite{benci_fortunato_1979}, \cite{berger_schechter_1972}, \cite{ASNSP_1973_3_27_4_591_0}, \cite{simader1996}, \cite{thater2002}, \cite{galdi_simader_1990}, and \cite{meyers_serrin_1960}. Alternatively, we can use homogeneous Sobolev spaces. \par 
In a paper published in 2019 \cite{leoni_tice_2019}, the trace space for the homogeneous Sobolev space $\dot W^{m,p}(\Omega)$ was determined for domains of type (\ref{generalomega}). A screened homogeneous fractional Sobolev space was introduced to capture the trace space, equipped with the seminorm, for $\Omega\subseteq \mathbb R^{N-1}$ open,
$$
|g|_{\tilde{W}^{s,p}_{\leq\eta}(\Omega)}:=\left(\int_{\Omega}\int_{\Omega\cap B(x,\eta(x))} \frac{|g(x)-g(y)|^p}{\Vert x-y\Vert^{N+sp}}dydx\right)^{1/p}$$
We note that in \cite{leoni_tice_2019}, they use the notation $|g|_{\tilde{W}^{s,p}_{(\eta)}(\Omega)}$; the change in notation is to make a clearer distinction with a different type of screened space that we discuss later. The core theorem proven in that paper was the following.
\begin{theorem}[Leoni, Tice \cite{leoni_tice_2019}]\label{leoniThm}
Let $\Omega$ be defined in (\ref{generalomega}) where $\eta^\pm:\mathbb R^{N-1}\rightarrow \mathbb R$ are Lipschitz and $\eta^-<\eta^+$. Let $L=|\eta^-|_{0,1}+|\eta^+|_{0,1}$, and let $1<p<\infty$.\begin{enumerate}
    \item 
 Then there exists $c=c(N,p)>0$ such that, for all $u\in \dot W^{1,p}(\Omega)$,
 \begin{align*}
\int_{\mathbb R^{N-1}}  \frac{\left\vert \Tr(u)(x',\eta^+(x'))-\Tr(u)(x',\eta^-(x'))\right\vert^p}{\left(\eta^+(x')-\eta^-(x')\right)^{sp-1}}\;dx'&\leq c\int_{\Omega}\left\vert \partial_Nu(x)\right\vert^pdx,\\
\left\vert \Tr(u)(\cdot,\eta^\pm(\cdot))\right\vert_{\tilde{W}^{1-1/p,p}_{\leq(\eta^+-\eta^-)/(2L)}(\mathbb R^{N-1})}^p&\leq (1+L)^p c\int_\Omega |\nabla u(x)|^pdx
 \end{align*}
\item Let $f^\pm\in L^p_{\text{Loc}}(\mathbb R^{N-1})$ be such that
\begin{align*}
\int_{\mathbb R^{N-1}}  \frac{\left\vert f^+(x')-f^-(x')\right\vert^p}{\left(\eta^+(x')-\eta^-(x')\right)^{sp-1}}\;dx'<\infty,\\
\left\vert f^{\pm}\right\vert_{ \tilde{W}^{1-1/p,p}_{\leq(\eta^+-\eta^-)/(2L)}(\mathbb R^{N-1})}<\infty.
\end{align*}
Then there exists $u\in \dot W^{1,p}(\Omega)$ such that $\Tr(u)(x',\eta^\pm)=f^\pm(x')$. Furthermore, the map $(f^+,f^-)\rightarrow u$ is linear, and, for $c=c(N,p)>0$,
\begin{align*}
\firstline|u|_{\dot W^{s,p}(\Omega)}^p\leq c(1+L)^p\int_{\mathbb R^{N-1}}  \frac{\left\vert f^+(x')-f^-(x')\right\vert^p}{\left(\eta^+(x')-\eta^-(x')\right)^{sp-1}}\;dx'\\
&\tab\tab+c(1+L)^p\left\vert f^{+}\right\vert_{ \tilde{W}^{1-1/p,p}_{\leq(\eta^+-\eta^-)/(2L)}(\mathbb R^{N-1})}^p+c(1+L)^p\left\vert f^{-}\right\vert_{ \tilde{W}^{1-1/p,p}_{\leq(\eta^+-\eta^-)/(2L)}(\mathbb R^{N-1})}^p
\end{align*} 
\end{enumerate}
\end{theorem}
Also see \cite{stevenson_tice_2020} for more information on screened Sobolev spaces. Historically, the trace space of fractional Sobolev spaces have been similar to their Sobolev space counterparts. Recall theorem \ref{generalset}, which showed $\Tr(W^{s,p}(\Omega))=W^{s-1/p,p}(\partial \Omega)$ for $\partial\Omega$ uniformly Lipschitz. A similar theorem states that $\Tr(W^{1,p}(\Omega))=W^{1-1/p,p}(\partial \Omega)$, showing that the trace space is the same. Therefore, the expectation would be that, similar to theorem \ref{leoniThm}, the trace space for $\dot W^{s,p}(\Omega)$ as in (\ref{generalomega}) is $\tilde{W}^{s-1/p,p}_{\leq (\eta^+-\eta^-)/2L}(\R^{N-1})$. However, we surprisingly find that this is not the case, and that we need to introduce a novel seminorm, which we will refer to as the \textbf{far screened Sobolev seminorm}, defined as
$$
|g|_{\tilde{W}^{s,p}_{\geq \eta}(\Omega)}:=\left(\int_{\Omega}\int_{\Omega\setminus B(x,\eta(x))} \eta(x')\eta(y')\frac{|g(x)-g(y)|^p}{\Vert x-y\Vert^{N+sp}}dydx\right)^{1/p},$$
and we define $\tilde{W}^{s,p}_{\geq \eta}(\Omega)$ to be the space for which the above seminorm is finite. Our trace results can now be stated. We first discuss the simpler case for $\eta^\pm$ constant. In other words, the set
\begin{equation}\label{flatomega}
\Omega:=\{(x',x_N)|b^+<x_N<b^-\}.
\end{equation}
Over this set, we derive the following theorem.
\begin{theorem}\label{theorem:flatMT}
Let $\Omega$ be defined in (\ref{flatomega}). Let $0<s<1$, $1<p<\infty$ and $1<sp$.
\begin{enumerate}
    \item 
 Then there exists $C=C(N,s,p,b^+-b^-)>0$ such that, for all $u\in \dot W^{s,p}(\Omega)$,
 \begin{align*}
\int_{\mathbb R^{N-1}}  \left\vert\Tr(u)(x',b^+)-\Tr(u)(x',b^-)\right\vert^p\;dx'\leq &C\left\vert u\right\vert_{\dot W^{s,p}(\Omega)}^p,\\
\left\vert \Tr(u)(\cdot,b^\pm)\right\vert_{\tilde{W}^{s-1/p,p}_{\leq b^+-b^-}(\mathbb R^{N-1})}\leq &C\left\vert u\right\vert_{\dot W^{s,p}(\Omega)},\\
\left\vert \Tr(u)(\cdot,b^\pm)\right\vert_{\tilde{W}^{s+1/p,p}_{\geq b^+-b^-}(\mathbb R^{N-1})}\leq &C\left\vert u\right\vert_{\dot W^{s,p}(\Omega)}
 \end{align*}
\item Let $f^\pm\in L^p_{\text{Loc}}(\mathbb R^{N-1})$ be such that
 \begin{align*}
\int_{\mathbb R^{N-1}}  \left\vert f^+(x')-f^-(x')\right\vert^p\;dx'< &\infty,\\
\left\vert f^\pm \right\vert_{\tilde{W}^{s-1/p,p}_{\leq b^+-b^-}(\mathbb R^{N-1})}< &\infty,\\
\left\vert f^\pm\right\vert_{\tilde{W}^{s+1/p,p}_{\geq b^+-b^-}(\mathbb R^{N-1})}< &\infty
 \end{align*}
Then there exists $u\in \dot W^{s,p}(\Omega)$ such that $\Tr(u)(x',b^\pm)=f^\pm(x')$. Furthermore, the map $(f^+,f^-)\rightarrow u$ is linear, and, for $C=C(N,s,p,b)>0$,
\begin{align*}
\firstline|u|_{\dot W^{s,p}(\Omega)}\leq C\Vert f^+-f^-\Vert_{L^p(\mathbb R^{N-1})}+C\left\vert f^+ \right\vert_{\tilde{W}^{s-1/p,p}_{\leq b^+-b^-}(\mathbb R^{N-1})}+C\left\vert f^- \right\vert_{\tilde{W}^{s-1/p,p}_{\leq b^+-b^-}(\mathbb R^{N-1})}\\
&\tab+C\left\vert f^+ \right\vert_{\tilde{W}^{s+1/p,p}_{\geq b^+-b^-}(\mathbb R^{N-1})}+C\left\vert f^- \right\vert_{\tilde{W}^{s+1/p,p}_{\geq b^+-b^-}(\mathbb R^{N-1})}
\end{align*}
\end{enumerate}
\end{theorem}
The principal strategy for addressing both parts of this simplified theorem is to define a seminorm equivalent to $|\cdot |_{\dot W^{s,p}(\Omega)}$. These equivalent seminorms look at slicing in specific dimensions, not unlike using Tonelli's theorem to separate the dimensions of an integral. This is especially useful in this simplified scenario, as the set $\{x_i:(x',x_i)\in \Omega\}\subset \mathbb R$ is always an interval. To address the second part, we mollify our function to generate a $C^{\infty}(\Omega)$ function that satisfies our desired properties. 
\par In Section $5$, we prove that, for certain parameters $s$, $p$, and $N$, we have the strict inclusions
$$\dot W^{s-1/p,p}(\R^N)\subsetneqq \tilde{W}^{s-1/p,p}_{\leq b^+-b^-}(\R^N)\cap \tilde{W}^{s+1/p,p}_{\geq b^+-b^-}(\R^N)\subsetneqq \tilde{W}^{s-1/p,p}_{\leq b^+-b^-}(\R^N).$$
Thus, the trace space of $\dot W^{s,p}(\Omega)$ does not coincide in general with the screened Sobolev space discussed in \cite{leoni_tice_2019}.\par 
After, we discuss the general theorem for $\Omega$ as in (\ref{generalomega}) and adapt similar techniques. 
\begin{theorem}\label{theorem:generalMT}
Let $\Omega$ be defined in (\ref{generalomega}) where $\eta^\pm:\mathbb R^{N-1}\rightarrow \mathbb R$ are Lipschitz and $\eta^-<\eta^+$. Let $L=|\eta^-|_{0,1}+|\eta^+|_{0,1}$. Let $1<p<\infty$ and $0<s<1$ such that $1<sp$.\begin{enumerate}
    \item 
 Then there exists $C=C(N,s,p,L)>0$ such that, for all $u\in \dot W^{s,p}(\Omega)$,
 \begin{align*}
\int_{\mathbb R^{N-1}}  \frac{\left\vert \Tr(u)(x',\eta^+(x'))-\Tr(u)(x',\eta^-(x'))\right\vert^p}{\left(\eta^+(x')-\eta^-(x')\right)^{sp-1}}\;dx'&\leq C\left\vert u\right\vert_{\dot W^{s,p}(\Omega)}^p,\\
\left\vert \Tr(u)(\cdot,\eta^\pm(\cdot))\right\vert_{\tilde{W}^{s-1/p,p}_{\leq(\eta^+-\eta^-)/(2L)}(\mathbb R^{N-1})}&\leq C\left\vert u\right\vert_{\dot W^{s,p}(\Omega)},\\
\left\vert \Tr(u)(\cdot,\eta^\pm(\cdot))\right\vert_{\tilde{W}^{s+1/p,p}_{\geq(\eta^+-\eta^-)/(2L)}(\mathbb R^{N-1})}&\leq C\left\vert u\right\vert_{\dot W^{s,p}(\Omega)},\\
 \end{align*}
\item Let $f^\pm\in L^p_{\text{Loc}}(\mathbb R^{N-1})$ be such that
\begin{align*}
\int_{\mathbb R^{N-1}}  \frac{\left\vert f^+(x')-f^-(x')\right\vert^p}{\left(\eta^+(x')-\eta^-(x')\right)^{sp-1}}\;dx'<\infty,\\
\left\vert f^{\pm}\right\vert_{ \tilde{W}^{s-1/p,p}_{\leq(\eta^+-\eta^-)/(2L)}(\mathbb R^{N-1})}<\infty,\\
\left\vert f^{\pm}\right\vert_{ \tilde{W}^{s+1/p,p}_{\geq(\eta^+-\eta^-)/(2L)}(\mathbb R^{N-1})}<\infty
\end{align*}
Then there exists $u\in \dot W^{1,p}(\Omega)$ such that $\Tr(u)(x',\eta^\pm)=f^\pm(x')$. Furthermore, the map $(f^+,f^-)\rightarrow u$ is linear, and, for $C=C(N,s,p,L)$,
\begin{align*}
\firstline|u|_{\dot W^{s,p}(\Omega)}^p\leq C\int_{\mathbb R^{N-1}}  \frac{\left\vert f^+(x')-f^-(x')\right\vert^p}{\left(\eta^+(x')-\eta^-(x')\right)^{sp-1}}\;dx'\\
&\tab\tab+C\left\vert f^{+}\right\vert_{ \tilde{W}^{s-1/p,p}_{\leq(\eta^+-\eta^-)/(2L)}(\mathbb R^{N-1})}^p+C\left\vert f^{-}\right\vert_{ \tilde{W}^{s-1/p,p}_{\leq(\eta^+-\eta^-)/(2L)}(\mathbb R^{N-1})}^p\\
&\tab\tab+C\left\vert f^{+}\right\vert_{ \tilde{W}^{s+1/p,p}_{\geq(\eta^+-\eta^-)/(2L)}(\mathbb R^{N-1})}^p+C\left\vert f^{-}\right\vert_{ \tilde{W}^{s+1/p,p}_{\geq(\eta^+-\eta^-)/(2L)}(\mathbb R^{N-1})}^p
\end{align*} 
\end{enumerate}
\end{theorem}
\par This paper is organized as follows. In Section $2$, we introduce relevant spaces, their intrinsic seminorms, and preliminary properties. In Section $3$, we prove theorem \ref{theorem:flatMT}. In Section $4$, we prove theorem \ref{theorem:generalMT}. In Section $5$, we study basic properties of the trace space and their relationship to previously established screened Sobolev spaces.

%% file: sections/2.Preliminaries.tex
\section{Preliminaries}

\subfile{2.Preliminaries/2.1.Notation}

\subfile{2.Preliminaries/2.2.EquivalentDomains}

\subfile{2.Preliminaries/2.3.SlicingLemmas}

\subfile{2.Preliminaries/2.4.FlatOmega}

\subfile{2.Preliminaries/2.5.GeneralOmega}

\subfile{2.Preliminaries/2.6.MollifierBounds}

%% file: sections/2.Preliminaries/2.1.Notation.tex
\subsection{Notation} \hfill\\
Let us first clarify our notation. First, we recall the definitions in our introduction. For $\Omega\subseteq \R^N$ open, $0<s<1$, $1\leq p<\infty$, we say $u\in \dot W^{s,p}(\Omega)$, the homogeneous fractional Sobolev space, if and only if
$$|u|_{\dot W^{s,p}(\Omega)}=\left(\int_\Omega \int_\Omega \frac{|u(x)-u(y)|^p}{\Vert x-y\Vert^{N+sp}}dydx\right)^{1/p}<\infty.$$
For $\eta$ Lipschitz, we say $u\in \dot W^{s,p}_{\leq \eta}(\Omega)$, the close screened homogeneous fractional Sobolev space, if and only if
$$|g|_{\tilde{W}^{s,p}_{\leq\eta}(\Omega)}:=\left(\int_{\Omega}\int_{\Omega\cap B(x,\eta(x))} \frac{|g(x)-g(y)|^p}{\Vert x-y\Vert^{N+sp}}dydx\right)^{1/p}<\infty.$$
We say $u\in \bar W^{s,p}_{\geq \eta}(\Omega)$, the far screened homogeneous fractional Sobolev space, if and only if
$$|g|_{\tilde{W}^{s,p}_{\geq \eta}(\Omega)}:=\left(\int_{\Omega}\int_{\Omega\setminus B(x,\eta(x))} \eta(x)\eta(y)\frac{|g(x)-g(y)|^p}{\Vert x-y\Vert^{N+sp}}dydx\right)^{1/p}<\infty.$$
Therefore, we have three different Sobolev-type seminorms. \newline 
We define
$$\Delta_{h}u(x)=u(x+h)-u(x).$$
\newline In order to discuss the H\"{o}lder continuous bound, for $f:\Omega\rightarrow \mathbb R$, we use
$$|f|_{0,\alpha}:=\sup\left\{\frac{|f(x)-f(y)|}{\Vert x-y\Vert^\alpha}:x,y\in \Omega\right\}.$$
Note that if $\alpha=1$, then this is just the Lipschitz bound.

%% file: sections/2.Preliminaries/2.2.EquivalentDomains.tex
\subsection{Equivalent Domains} \hfill\\
We first desire to simplify the domains that we hope to work with. Namely, we namely wish to work with the domain
\begin{equation}\label{generalomegawflat}
\Omega=\{(x',x_N)|0<x_N<\eta(x')\}.
\end{equation}
for $\eta$ $1$-Lipschitz. Therefore, we start with some lemmas that map sets in (\ref{generalomega}) to (\ref{generalomegawflat}). First, we want to be able to flatten the boundary, which leads us to the first theorem allowing translations in the $x_N$ dimension.
\begin{lemma}\label{lem:translation}
Let $\Omega$ open, and $\eta(x')$ Lipschitz continuous such that $|\eta|_{0,1}=L$. We define
$$\Omega':=\{(x',x_N)|(x',x_N-\eta(x_N))\in \Omega\}.$$
Then $u\in \dot W^{s,p}(\Omega)$ if and only if $v\in \dot W^{s,p}(\Omega')$ for $v$ defined as
$$v(x',x_N)=u(x',x_N+\eta(x')).$$
Furthermore, there exists $C_1,C_2>0$ dependent on $s$, $p$, and $L$ such that
$$C_1|v|_{\dot W^{s,p}(\Omega')}\leq |u|_{\dot W^{s,p}(\Omega)}\leq C_2|v|_{\dot W^{s,p}(\Omega')}.$$
\end{lemma}
\begin{proof}
First, we show $v\in \dot W^{s,p}(\Omega')$. We consider:
\begin{align*}
|v|_{\dot W^{s,p}(\Omega')}^p&=\int_{\Omega'}\int_{\Omega'} \frac{|v(x)-v(y)|^p}{\Vert x-y\Vert^{N+sp}}dydx\\
&=\int_{\Omega'}\int_{\Omega'} \frac{|u(x',x_N+\eta(x'))-u(y',y_N+\eta(y'))|^p}{\Vert x-y\Vert^{N+sp}}dydx.
\end{align*}
We note
\begin{align*}
\Vert (x',x_N+\eta(x'))-(y',y_N+\eta(y'))\Vert&\leq \Vert x'-y'\Vert +|x_N-y_N|+|\eta(x')-\eta(y')|\\
&\leq \Vert x'-y'\Vert + |x_N-y_N|+L\Vert x'-y'\Vert\\
&\leq (2+L)\Vert x-y\Vert.
\end{align*}
Therefore, we get
$$|v|_{\dot W^{s,p}(\Omega')}^p\geq C\int_{\Omega'}\int_{\Omega'} \frac{|u(x',x_N+\eta(x'))-u(y',y_N+\eta(y'))|^p}{\Vert (x',x_N+\eta(x'))-(y',y_N+\eta(y'))\Vert^{N+sp}}dydx.$$
With a change of variables $x_N=x_N+\eta(x')$ and $y_N=y_N+\eta(y')$, we get:
$$=C\int_{\Omega}\int_{\Omega} \frac{|u(x)-u(x)|^p}{\Vert x-y\Vert^{N+sp}}dydx.$$
This proves the lower bound; since this proof only relied on $\eta(x')$ being Lipschitz, if we applied the computations with $-\eta(x')$ starting with $u$ instead, we would instead get:
$$|u|_{\dot W^{s,p}(\Omega)}^p\geq C\int_{\Omega'}\int_{\Omega'} \frac{|v(x)-v(x)|^p}{\Vert x-y\Vert^{N+sp}}dydx.$$
\end{proof}
Therefore, note that for $\Omega$ as in (\ref{generalomega}), we can instead consider
$$\Omega':=\{(x',x_N):0<x_N<\eta(x')\}$$
by applying lemma~\ref{lem:translation}, for $\eta(x'):=\eta^+(x')-\eta^-(x')$. Next, we wish to control the Lipschitz bound on $\eta$. We note that for $\eta$ Lipschitz and non-constant, the function $\frac{\eta}{|\eta|_{0,1}}$ is $1$-Lipschitz. This leads us to the following theorem for scaling domains.
\begin{lemma}\label{lem:dilation}
Let $\Omega$ open, and let $\alpha>0$. We define
$$\Omega':=\{(x',x_N)|(x',\alpha x_N)\in \Omega\}.$$
Then $u\in \dot W^{s,p}(\Omega)$ if and only if $v\in \dot W^{s,p}(\Omega')$ for $v$ defined as
$$v(x',x_N)=u(x',\alpha x_N).$$
Furthermore, there exists $C_1,C_2>0$ dependent on $s$, $p$, and $\alpha$ such that
$$C_1|v|_{\dot W^{s,p}(\Omega')}\leq |u|_{\dot W^{s,p}(\Omega)}\leq C_2|v|_{\dot W^{s,p}(\Omega')}.$$
\end{lemma}
\begin{proof}
We would like to compare $\Vert x-y\Vert$ with $\Vert (x',\alpha x_N)-(y',\alpha y_N)\Vert$. We observe
\begin{align*} 
\Vert (x',\alpha x_N)-(y',\alpha y_N)\Vert& \leq \Vert x'-y'\Vert+\alpha |x_N-y_N|\\
&\leq (1+\alpha) \Vert x-y\Vert.
\end{align*}
Therefore,
\begin{align*}
|v|_{\dot W^{s,p}(\Omega')}^p&=\int_{\Omega'}\int_{\Omega'} \frac{|v(x)-v(y)|^p}{\Vert x-y\Vert^{N+sp}}dydx\\
&=\int_{\Omega'}\int_{\Omega'} \frac{|u(x',\alpha x_N)-v(y',\alpha y_N)|^p}{\Vert x-y\Vert^{N+sp}}dydx\\
&\geq C\int_{\Omega'}\int_{\Omega'} \frac{|u(x',\alpha x_N)-v(y',\alpha y_N)|^p}{\Vert (x',\alpha x_N)-(y',\alpha y_N)\Vert^{N+sp}}dydx.
\end{align*}
We use the change of variables $x_N=\alpha x_N$ and $y_N=\alpha y_N$.
\begin{align*}
&= C\int_{\Omega}\int_{\Omega} \frac{|u(x',x_N)-v(y',y_N)|^p}{\Vert x-y\Vert^{N+sp}}dydx.
\end{align*}
We again use $1/\alpha$ to derive the other bound. 
\end{proof}
Therefore, we can always scale by $1/|\eta|_{0,1}$ in order to get $1$-Lipschitz function. Lastly, we consider a natural property of inclusion. Note that if $\Omega'\subset \Omega$ where $\Omega'$ and $\Omega$ open, then naturally, $f\in \dot W^{s,p}(\Omega)$ implies $f\restriction \Omega'\in \dot W^{s,p}(\Omega')$. Therefore, if we were to consider the sets:
$$\Omega_1:=\{(x',x_N)|0<x_N<\eta_1(x')\}$$
$$\Omega_2:=\{(x',x_N)|0<x_N<\eta_2(x')\}.$$
such that $0<\eta_2<\eta_1$, any bound we impose on $\Tr(f)(\cdot,0)$ for $f\in \dot W^{s,p}(\Omega)$, then that restriction is still true for $\Tr(f\restriction \Omega_2)(\cdot,0)$ if we treat $f\restriction \Omega_2\in W^{s,p}(\Omega_2)$. Therefore, we first prove a property derived from this intuition. 
\begin{lemma}\label{lem:decreaser}
For $0<\eta_1$ Lipschitz, $0<\eta_2<\eta_1$ Lipschitz with constant $L$, if $g$ measurable such that:
$$|g|_{\tilde{W}^{s-1/p,p}_{\leq \eta_1}(\mathbb R^{N-1})}<\infty$$
$$|g|_{\tilde{W}^{s+1/p,p}_{\geq \eta_1}(\mathbb R^{N-1})}<\infty$$
Then we have:
$$|g|_{\tilde{W}^{s-1/p,p}_{\leq \eta_2}(\mathbb R^{N-1})}<\infty$$
$$|g|_{\tilde{W}^{s+1/p,p}_{\geq \eta_2}(\mathbb R^{N-1})}<\infty$$
\end{lemma}
\begin{proof}
Clearly:
$$|g|_{\tilde{W}^{s-1/p,p}_{\leq \eta_2}(\mathbb R^{N-1})}\leq |g|_{\tilde{W}^{s-1/p,p}_{\leq \eta_1}(\mathbb R^{N-1})}<\infty$$
Therefore, we just need to consider the other term. We have
\begin{align*}
|g|_{\tilde{W}^{s+1/p,p}_{\geq \eta_2}(\mathbb R^{N-1})}^p&=\int_{\mathbb R^{N-1}}\int_{\mathbb R^{N-1}\setminus B(x',\eta_2(x'))} \eta_2(x')\eta_2(y')\frac{|g(x')-g(y')|^p}{\Vert x'-y'\Vert^{N+sp}}dy'dx'\\
&\leq \int_{\mathbb R^{N-1}}\int_{\mathbb R^{N-1}\setminus B(x',\eta_1(x'))} \eta_1(x')\eta_1(y')\frac{|g(x')-g(y')|^p}{\Vert x'-y'\Vert^{N+sp}}dy'dx'\\
&\tab+\int_{\mathbb R^{N-1}}\int_{B(x',\eta_1(x'))\setminus B(x',\eta_2(x'))} \eta_2(x')\eta_2(y')\frac{|g(x')-g(y')|^p}{\Vert x'-y'\Vert^{N+sp}}dy'dx'.
\end{align*}
The first term is $|g|_{\tilde{W}^{s+1/p,p}_{\geq \eta_1}(\mathbb R^{N-1})}^p$.  Meanwhile, for the second term, we note that, using the set we are integrating over: 
$$\eta_2(x')\leq \Vert x'-y'\Vert.$$ 
On the other hand, using the Lipschitz bound:
$$\eta_2(y')\leq \eta_2(x')+L\Vert x'-y'\Vert \leq \left(1+L\right) \Vert x'-y'\Vert.$$
Therefore, the second term is bounded by $(L+1)|g|_{\tilde{W}^{s-1/p,p}_{\leq \eta_1}(\mathbb R^{N-1})}^p$ for some $C$. 
\end{proof}

%% file: sections/2.Preliminaries/2.3.SlicingLemmas.tex
\subsection{Slicing Lemmas} \hfill\\ 
We first develop tools to traverse up or down dimensions while working with fractional Sobolev spaces. This is mainly important in separating the dimensions parallel to the sides of the strip, and the dimension perpendicular. 
\begin{lemma}\label{lem:slicingplane}
Given $\R^N$ and $\lambda>N$, there exists constant $C_1=C_1(\lambda,N)$ and $C_2=C_2(\lambda,N)$ such that:
$$\frac{C_1}{|x_N-y_N|^{\lambda-N+1}}\leq \int_{\mathbb R^{N-1}} \frac{1}{\Vert x-y\Vert^\lambda} dx'_N\leq  \frac{C_1}{|x_N-y_N|^{\lambda-N+1}}.$$
\end{lemma}
\begin{proof}
First, consider, for a fixed $y\in \R$:
\begin{align*}
\int_\R \frac{1}{\rho^\lambda + |x-y|^\lambda} dx
&=2\int_0^\infty \frac{1}{\rho^\lambda+t^\lambda} dt\\
&=\frac{2}{\rho^{\lambda-1}}\int_0^\infty \frac{1}{1+t^\lambda} dt.
\end{align*}
Therefore, we get some $C=C(\lambda)$ such that:
$$\int_\R \frac{1}{\rho^\lambda+|x-y|^\lambda} dx = \frac{C}{\rho^{\lambda-1}}.$$
Now, reconsider the main problem. We will continue with induction on $N$. If $N=2$, we have that there are $C_1,C_2$ dependent on $\lambda$ such that:
$$\frac{C_1}{|x_1-y_1|^\lambda + |x_2-y_2|^\lambda}\leq \frac{1}{\Vert x-y\Vert^\lambda}\leq \frac{C_2}{|x_1-y_1|^\lambda + |x_2-y_2|^\lambda}.$$
By integrating over $\mathbb R$ and then applying the previous part, we get:
$$\frac{C_1\cdot C}{|x_2-y_2|^{\lambda-1}}\leq \int_{\mathbb R} \frac{1}{\Vert x-y\Vert^\lambda} dx_1\leq  \frac{C_2\cdot C}{|x_2-y_2|^{\lambda-1}}.$$
Now, assume the result holds for $N-1$. Again, note that we can find $C_1,C_2$ dependent on $N$ such that:
$$\frac{C_1}{\Vert x'_1-y'_1\Vert^\lambda + |x_1-y_1|^\lambda}\leq \frac{1}{\Vert x-y\Vert^\lambda}\leq \frac{C_2}{\Vert x'_1-y'_1\Vert^\lambda + |x_1-y_1|^\lambda}.$$
Therefore, integrating and using $\rho=\Vert x'_1-y'_1\Vert$ in the first part, we get:
$$\frac{C_1\cdot C}{\Vert x'-y'_1\Vert^{\lambda-1}}\leq \int_{\mathbb R} \frac{1}{\Vert x-y\Vert^\lambda} dx_1\leq  \frac{C_2\cdot C}{\Vert x'_1-y'_1\Vert^{\lambda-1}}.$$
Now integrate both sides over $x_2,...,x_{N-1}$ and then applying the inductive hypothesis, we get the desired result. 
\end{proof}
Next, we look at when the intervals are bounded. 
\begin{lemma}\label{lem:slicingfinite}
Let $r>0$, $k>0$, $\lambda>1$, $0<\rho\leq kr$. There exists $C_1=C_1(\lambda,k)>0$ and $C_2=C_2(\lambda,k)>0$ such that for every $y\in (0,r)$:
    $$\frac{C_1}{\rho^{\lambda-1}}\leq \int_0^r \frac{1}{\rho^{\lambda}+|x-y|^\lambda}dx\leq \frac{C_2}{\rho^{\lambda-1}}.$$
\end{lemma}
\begin{proof}We have:
\begin{align*}
\int_0^r \frac{1}{\rho^\lambda+|x-y|^\lambda}dx
&=\int_0^y \frac{1}{\rho^\lambda+t^\lambda}dt+\int_0^{r-y} \frac{1}{\rho^\lambda+t^\lambda}dt\\
&= \frac{1}{\rho^{\lambda-1}} \int_0^{y/\rho} \frac{1}{1+t^\lambda} dt+\frac{1}{\rho^{\lambda-1}}\int_0^{(r-y)/\rho} \frac{1}{1+t^\lambda} dt.
\end{align*}
Now, depending on $y$ with respects to $r/2$, we have that one of the two integrals is at least:
$$\frac{1}{\rho^{\lambda-1}} \int_0^{r/2\rho} \frac{1}{1+t^\lambda} dt.$$
Therefore, we get:
$$\frac{1}{\rho^{\lambda-1}}\int_0^{r/2\rho} \frac{1}{1+t^\lambda} dt \leq \int_0^r \frac{1}{\rho^\lambda+|x-y|^\lambda} dx\leq \frac{2}{\rho^{\lambda-1}} \int_0^\infty \frac{1}{1+t^\lambda} dt.$$
Since $\rho<kr$, we can simplify the lower bound to:
$$\frac{1}{\rho^{\lambda-1}}\int_0^{1/2k} \frac{1}{1+t^\lambda} dt \leq \int_0^r \frac{1}{\rho^\lambda+|x-y|^\lambda} dx\leq \frac{2}{\rho^{\lambda-1}} \int_0^\infty \frac{1}{1+t^\lambda} dt.$$
\end{proof}

%% file: sections/2.Preliminaries/2.4.FlatOmega.tex
\subsection{Equivalent Seminorms for Strip Domains} \hfill\\
We wish to define seminorms that are equivalent to the standard homogeneous fractional Sobolev space seminorms. Similar to how we use Tonelli to reduce an integral to integrating over dimensions, we can reduce certain seminorms to integrating only in certain dimensions. We first consider $\Omega$ as in (\ref{flatomega}). 
\begin{theorem}\label{Thm:BoundOnSliceFlat}
Let $\Omega$ be as in (\ref{flatomega}). Then there is $C=C(N,s,p)>0$ such that for all $u\in \dot W^{s,p}(\Omega)$:
\begin{enumerate}
    \item \[\int_{\R^{N-1}} |u(x',\cdot)|_{\dot W^{s,p}((0,b))}^p dx'<C|u|_{\dot W^{s,p}(\Omega)}^p\]
    \item \[\int_0^b |u(\cdot,x_N)|_{\dot W^{s,p}_{\leq b}(\mathbb R^{N-1})}^pdx_N<C|u|_{\dot W^{s,p}(\Omega)}^p\]
    \item \[\int_0^b |u(\cdot,x_N)|_{\dot W^{s+1/p,p}_{\geq b}(\mathbb R^{N-1})}^pdx_N<Cb|u|_{\dot W^{s,p}(\Omega)}^p\]
\end{enumerate}
\end{theorem}
\begin{proof}
We look at each individual seminorm separately. \hfill\\
\textbf{Part 1:} Let $Q(x',x_N,y_N)$ denote the cube with center $\left(x',\frac{x_N+y_N}{2}\right)$ and side length $\frac{1}{4\sqrt{N}}|x_N-y_N|$. Then, we have
\begin{align*}
\firstline\int_{\R^{N-1}} |u(x',\cdot)|_{\dot W^{s,p}((0,b))}^pdx'\\
&=\int_{\R^{N-1}}\int_0^b\int_0^b \frac{|u(x',x_N)-u(x',y_N)|^p}{|x_N-y_N|^{1+sp}}dy_Ndx_Ndx'\\
&\leq C\int_{\R^{N-1}}\int_0^b\int_0^b\frac{1}{|x_N-y_N|^N}\int_{Q(x',x_N,y_N)} \frac{|u(x',x_N)-u(x',y_N)|^p}{|x_N-y_N|^{1+sp}}dzdy_Ndx_Ndx'\\
&\leq C\int_{\R^{N-1}}\int_0^b\int_0^b\int_{Q(x',x_N,y_N)} \frac{|u(x)-u(z)|^p}{|x_N-y_N|^{N+1+sp}}dzdy_Ndx_Ndx'\\
&=C\int_\Omega\int_\Omega |u(x)-u(z)|^p \int_0^b \chi_{Q(x',x_N,y)}(z)\frac{1}{|x_N-y_N|^{N+1+sp}} dy_N dz dx.
\end{align*}
Now that since $z\in Q(x',x_N,y_N)$, then: 
$$\left \Vert z-\left(x',\frac{x_N+y_N}{2}\right)\right\Vert \leq \frac{|x_N-y_N|}{4}.$$
Furthermore, we have:
$$\left\Vert x-\left(x',\frac{x_N+y_N}{2}\right)\right\Vert=\frac{|x_N-y_N|}{2}.$$
Therefore, with triangle inequality, we have:
\begin{equation}\label{cubebound}
\frac{4}{3}\Vert x-z\Vert \leq |x_N-y_N|.
\end{equation}
Thus, we can simplify to:
$$\leq C\int_\Omega\int_\Omega |u(x)-u(z)|^p \int_{\frac{4}{3}\Vert x-z\Vert}^\infty \frac{1}{|x_N-y_N|^{N+1+sp}} dy_N dz dx$$
$$\leq C\int_\Omega\int_\Omega \frac{|u(x)-u(z)|^p }{\Vert x-z\Vert^{N+sp}} dz dx.$$
\textbf{Part 2:}
Let $Q(x',y',x_N)$ denote the cube with corner at $\left ( \frac{x'+y'}{2},x_N\right)$ and side length $\frac{1}{4\sqrt{N}}\Vert x'-y'\Vert$. Since the side length is less than $b/2$, we can choose at least one direction for the cube to go such that the cube is entirely contained in $\Omega$. Then we get
\begin{align*}
\firstline \int_0^b |u(\cdot,x_N)|_{\dot W^{s,p}_{\leq b}(\mathbb R^{N-1})}^pdx_N\\
&\int_0^b \int_{\R^{N-1}}\int_{B(x',b)} \frac{|u(x',x_N)-u(y',x_N)|^p}{\Vert x'-y'\Vert^{N-1+sp}}dy'dx'dx_N\\
&\leq C\int_0^b\int_{\R^{N-1}}\int_{B(x',b)}\frac{1}{\Vert x'-y'\Vert^N}\int_{Q(x',y',x_N)} \frac{|u(x',x_N)-u(y',x_N)|^p}{\Vert x'-y'\Vert^{N-1+sp}}dzdy'dx'dx_N\\
&\leq C\int_0^b\int_{\R^{N-1}}\int_{B(x',b)}\int_{Q(x',y',x_N)} \frac{|u(x)-u(z)|^p}{\Vert x'-y'\Vert^{2N-1+sp}}dzdy'dx'dx_N\\
&=C\int_\Omega\int_\Omega |u(x)-u(z)|^p \int_{B(x',b)} \chi_{Q(x',y',x_N)}(z)\frac{1}{\Vert x'
-y'\Vert^{2N-1+sp}} dy' dz dx.
\end{align*}
By (\ref{cubebound}), we have 
$$\frac{4}{3}\Vert x-z\Vert \leq \Vert x'-y'\Vert.$$
Thus, we can simplify to
\begin{align*}
&\leq C\int_\Omega\int_\Omega |u(x)-u(z)|^p \int_{\R^{N-1}\setminus B(x',\Vert x-z\Vert)} \frac{1}{\Vert x'-y'\Vert^{2N-1+sp}} dy' dz dx\\
&\leq C\int_\Omega\int_\Omega \frac{|u(x)-u(z)|^p }{\Vert x-z\Vert^{N+sp}} dz dx.
\end{align*}
\textbf{Part 3:} Let $R(x',y',x_N)$ denote the rectangle with corner at $\left(\frac{x'+y'}{2},x_N\right)$, height $\frac{b}{4\sqrt{N}}$ along $(0,b)$, and width $\frac{1}{4\sqrt{N}}\Vert x'-y'\Vert$ along other directions. Therefore, we have
\begin{align*}
\firstline \int_0^b |u(\cdot,x_N)|_{\dot W^{s+1/p,p}_{\geq b}(\mathbb R^{N-1})}^pdx_N\\
&=\int_0^b \int_{\R^{N-1}}\int_{\R^{N-1}\setminus B(x',b)} \frac{|u(x',x_N)-u(y',x_N)|^p}{\Vert x'-y'\Vert^{N+sp}}dy'dx'dx_N\\
&\leq \frac{C}{b}\int_0^b\int_{\R^{N-1}}\int_{\R^{N-1}\setminus B(x',b)}\frac{1}{\Vert x'-y'\Vert^{N-1}}\int_{R(x',y',x_N)} \frac{|u(x',x_N)-u(y',x_N)|^p}{\Vert x'-y'\Vert^{N+sp}}dzdy'dx'dx_N\\
&\leq \frac{C}{b}\int_0^b\int_{\R^{N-1}}\int_{\R^{N-1}\setminus B(x',b)}\int_{R(x',y',x_N)} \frac{|u(x)-u(z)|^p}{\Vert x'-y'\Vert^{2N-1+sp}}dzdy'dx'dx_N\\
&=\frac{C}{b}\int_\Omega\int_\Omega |u(x)-u(z)|^p \int_{\R^{N-1}\setminus B(x',b)} \chi_{Q(x',y',x_N)}(z)\frac{1}{\Vert x'
-y'\Vert^{2N-1+sp}} dy' dz dx. 
\end{align*}
By (\ref{cubebound}), we have 
$$\frac{4}{3}\Vert x-z\Vert \leq \Vert x'-y'\Vert.$$
Thus, we can simplify to:
\begin{align*}
&\leq \frac{C}{b}\int_\Omega\int_\Omega |u(x)-u(z)|^p \int_{\R^{N-1}\setminus B(x',\Vert x-z\Vert)} \frac{1}{\Vert x'-y'\Vert^{2N-1+sp}} dy' dz dx\\
&\leq \frac{C}{b}\int_\Omega\int_\Omega \frac{|u(x)-u(z)|^p }{\Vert x-z\Vert^{N+sp}} dz dx.
\end{align*}
\end{proof}
\begin{theorem}\label{Thm:BoundonSobolevFlat}
Let $\Omega$ be as in (\ref{flatomega}). Let $u: \Omega\rightarrow \mathbb R$ such that
\begin{align*}
    \int_{\R^{N-1}} |u(x',\cdot)|_{\dot W^{s,p}((0,b))}^pdx'<&\infty,\\
    \int_0^b |u(\cdot,x_N)|_{\dot W^{s,p}_{\leq b}(\mathbb R^{N-1})}^pdx_N<&\infty,\\
    \int_0^b |u(\cdot,x_N)|_{\dot W^{s+1/p,p}_{\geq b}(\mathbb R^{N-1})}^pdx_N<&\infty.
\end{align*}
Then $u\in \dot W^{s,p}(\Omega)$. Furthermore, there exists $C=C(N,s,p)>0$ such that
\begin{align*}
C|u|_{\dot W^{s,p}(\Omega)}^p&\leq\int_{\R^{N-1}} |u(x',\cdot)|_{\dot W^{s,p}((0,b))}^p\;dx'\\
&\tab +\int_0^b |u(\cdot,x_N)|_{\dot W^{s,p}_{\leq b}(\mathbb R^{N-1})}Yp\;dx_N\\
&\tab + b\int_0^b |u(\cdot,x_N)|_{\dot W^{s+1/p,p}_{\geq b}(\mathbb R^{N-1})}^p\;dx_N.
\end{align*}
\end{theorem}
\begin{proof}
We have
\begin{align*}
|u|_{\dot W^{s,p}(\Omega)}^p&=\int_{\Omega}\int_{\Omega} \frac{|u(x)-u(y)|^p}{\Vert x-y\Vert^{N+sp}} dy dx\\
&=\int_{\R^{N-1}}\int_0^b \int_{\R^{N-1}}\int_0^b \frac{|u(x)-u(y)|^p}{\Vert x-y\Vert^{N+sp}} dy_N dy'dx_Ndx'\\ 
&\leq 2^{p-1} \int_{\R^{N-1}}\int_0^b \int_{\R^{N-1}}\int_0^b \frac{|u(x',x_N)-u(y',x_N)|^p}{\Vert x-y\Vert^{N+sp}} dy_N dy'dx_Ndx'\\
&\tab+ 2^{p-1} \int_{\R^{N-1}}\int_0^b \int_{\R^{N-1}}\int_0^b \frac{|u(y',x_N)-u(y',y_N)|^p}{\Vert x-y\Vert^{N+sp}} dy_N dy'dx_Ndx'\\
&=2^{p-1}\mathcal A+2^{p-1}\mathcal B.
\end{align*}
We consider each term separately. First consider
\begin{align*}
\mathcal A&=\int_{\R^{N-1}}\int_0^b \int_{B(x',b)}\int_0^b \frac{|u(x',x_N)-u(y',x_N)|^p}{\Vert x-y\Vert^{N+sp}} dy_N dy'dx_Ndx'\\
&\tab+\int_{\R^{N-1}}\int_0^b \int_{\R^{N-1}\setminus B(x',b)}\int_0^b \frac{|u(x',x_N)-u(y',x_N)|^p}{\Vert x-y\Vert^{N+sp}} dy_N dy'dx_Ndx'\\
&=\mathcal A_1+\mathcal A_2.
\end{align*}
Again, we will analyze each term separately. First, consider
$$\mathcal A_1=\int_{\R^{N-1}}\int_0^b \int_{B(x',b)} |u(x',x_N)-u(y',x_N)|^p \int_0^b \frac{1}{\Vert x-y\Vert^{N+sp}} dy_N dy'dx_Ndx'.$$
Note that by lemma \ref{lem:slicingfinite}, we get a bound of
$$C\int_{\R^{N-1}}\int_0^b \int_{B(x',b)} \frac{|u(x',x_N)-u(y',x_N)|^p}{\Vert x'-y'\Vert^{N-1+sp}} dy_N dy'dx_Ndx'.$$
Next, we consider
$$\mathcal A_2=\int_{\R^{N-1}}\int_0^b \int_{\R^{N-1}\setminus B(x',b)} |u(x',x_N)-u(y',x_N)|^p \int_0^b \frac{1}{\Vert x-y\Vert^{N+sp}} dy_N dy'dx_Ndx'.$$
Since $\Vert x'-y'\Vert\leq \Vert x-y\Vert$, we get
\begin{align*}
&\leq \int_{\R^{N-1}}\int_0^b \int_{\R^{N-1}\setminus B(x',b)} |u(x',x_N)-u(y',x_N)|^p \frac{b}{\Vert x'-y'\Vert^{N+sp}} dy_N dy'dx_Ndx'\\
&=b \int_{\R^{N-1}} \int_0^b \int_{\R^{N-1}\setminus B(x',b)} \frac{|u(x',x_N)-u(y',x_N)|^p}{\Vert x'-y'\Vert^{N+sp}} dy' dx_N dx'.
\end{align*}
Lastly, we consider
\begin{align*}
\mathcal B&=\int_{\R^{N-1}}\int_0^b \int_{\R^{N-1}}\int_0^b \frac{|u(y',x_N)-u(y',y_N)|^p}{\Vert x-y\Vert^{N+sp}} dy_N dy'dx_Ndx'\\
&=\int_{\mathbb R^{N-1}}\int_0^b\int_0^b |u(y',x_N)-u(y',y_N)|^p \int_{\mathbb R^{N-1}} \frac{1}{\Vert x-y\Vert^{N+sp}}\;dx'dx_Ndy_Ndy'.
\end{align*}
By lemma \ref{lem:slicingplane}:
$$\leq \int_{\mathbb R^{N-1}}\int_0^b\int_0^b \frac{|u(y',x_N)-u(y',y_N)|^p}{|x_N-y_N|^{1+sp}} dx_Ndy_Ndy'.$$
\end{proof}

%% file: sections/2.Preliminaries/2.5.GeneralOmega.tex
\subsection{Equivalent Seminorms for Strip-like Domains}\hfill\\
We hope to establish similar lemmas for $\Omega$ as in (\ref{generalomegawflat}); however, this poses several problems. Specifically, we note that a level set of $\Omega$ is no longer convex. Therefore, for example, we can no longer take rectangle subsets of $\Omega$ between arbitrary points. 
\begin{theorem}
Let $\Omega$ be as in (\ref{generalomegawflat}). There is $C=C(N,s,p)$ such that for all $u\in \dot W^{s,p}(\Omega)$,
\begin{enumerate}
    \item \[\int_{\R^{N-1}}\int_0^{\eta(x')}\int_0^{\eta(x')}\frac{|u(x',x_N)-u(x',y_N)|^p}{|x_N-y_N|^{1+sp}}dy_Ndx_Ndx'\leq C\vert u\vert_{\dot W^{s,p}(\Omega)}^p\]
    \item \[\int_{\R^{N-1}}\int_0^{\frac{\eta(x')}{2}}\int_{B(x',\frac{\eta(x')}{2} )} \frac{|u(x',x_N)-u(y',x_N)|^p}{\Vert x'-y'\Vert^{N-1+sp}}dy'dx_Ndx'\leq C \vert u\vert_{\dot W^{s,p}(\Omega)}^p\]
\end{enumerate}
\end{theorem}
\begin{proof}
\hfill\\
We consider the two parts separately.\hfill\\
\textbf{Part 1:} Let $Q(x',x_N,y_N)$ denote the cube with center $\left(x',\frac{x_N+y_N}{2}\right)$ and side length $\frac{1}{4\sqrt{N}}|x_N-y_N|$. We need to confirm that this cube is contained in $\Omega$. First, without loss of generality, assume $x_N>y_N$; we have that the set $\{(z',z_N)\in \mathbb R^N|0<z_N<x_N-\Vert x'-z'\Vert\}$ is in $\Omega$, using the Lipschitz bound. Now, note that for all $(z',z_N)\in Q(x',x_N,y_N)$, we have that:
$$z_N<\frac{1}{8\sqrt{N}}(x_N-y_N)+\frac{x_N+y_N}{2}.$$
Meanwhile, we have that:
$$\Vert x'-z'\Vert\leq \frac{1}{4}(x_N-y_N).$$
Therefore, we note that using these two inequalities, we do have that:
\begin{align*}
x_N-\Vert x'-z'\Vert &\geq x_N-\frac{1}{4}(x_N-y_N)\\
&=\frac{3}{4}x_N+\frac{1}{4}y_N\\
&=\frac{x_N+y_N}{2}+\frac{x_N-y_N}{4}>z_N.
\end{align*}
Since we have the cube is contained in $\Omega$, we continue. Then, we have:
\begin{align*}
\firstline \int_{\R^{N-1}}\int_0^{\eta(x')}\int_0^{\eta(x')} \frac{|u(x',x_N)-u(x',y_N)|^p}{|x_N-y_N|^{1+sp}}dy_Ndx_Ndx'\\
&\leq C\int_{\R^{N-1}}\int_0^{\eta(x')}\int_0^{\eta(x')}\frac{1}{|x_N-y_N|^N}\int_{Q(x',x_N,y_N)} \frac{|u(x',x_N)-u(x',y_N)|^p}{|x_N-y_N|^{1+sp}}dzdy_Ndx_Ndx'\\
&\leq C\int_{\R^{N-1}}\int_0^{\eta(x')}\int_0^{\eta(x')}\int_{Q(x',x_N,y_N)} \frac{|u(x)-u(z)|^p}{|x_N-y_N|^{N+1+sp}}dzdy_Ndx_Ndx'\\
&=C\int_\Omega\int_\Omega |u(x)-u(z)|^p \int_0^{\eta(x')} \chi_{Q(x',x_N,y)}(z)\frac{1}{|x_N-y_N|^{N+1+sp}} dy_N dz dx.
\end{align*}
By (\ref{cubebound}), we have 
\begin{align*}
\frac{4}{3}\Vert x-z\Vert \leq |x_N-y_N|. 
\end{align*}
Thus, we can simplify to:
\begin{align*}
&\leq C\int_\Omega\int_\Omega |u(x)-u(z)|^p \int_{\frac{4}{3}\Vert x-z\Vert}^\infty \frac{1}{|x_N-y_N|^{N+1+sp}} dy_N dz dx\\
&\leq C\int_\Omega\int_\Omega \frac{|u(x)-u(z)|^p }{\Vert x-z\Vert^{N+sp}} dz dx.
\end{align*}
\textbf{Part 2:} Let $Q(x',y',x_N)$ denote a cube that contains $\left ( \frac{x'+y'}{2},x_N\right)$ and side length $\frac{1}{4 \sqrt{N}}\Vert x'-y'\Vert$. We need to show that such a cube can always be found inside $\Omega$. First, note that if we consider $\Omega$, then the cylinder with base centered at $(x',0)$, radius $\frac{\eta(x')}{2}$, and height $\frac{\eta(x')}{2}$ is contained in $\Omega$. Note that for any point $z$ in such a cylinder, we can find a cube of side length $\frac{\eta(x')}{4}$ that contains $z$, and is fully contained in the cylinder. Note that our point $\left(\frac{x'+y'}{2},x_N\right)$ is contained in this cylinder, and we are looking for a cube with shorter side length; therefore, it is always possible. Specifically, note that we are not saying that such a cube needs to be axis-aligned, just that it needs to exist. Then, we have:
\begin{align*}
\firstline\int_{\R^{N-1}}\int_0^{\frac{\eta(x')}{2}}\int_{B(x',\frac{\eta(x')}{2} )} \frac{|u(x',x_N)-u(y',x_N)|^p}{\Vert x'-y'\Vert^{N-1+sp}}dy'dx'dx_N\\
&\leq C\int_{\R^{N-1}}\int_0^{\frac{\eta(x')}{2}}\int_{B(x',\frac{\eta(x')}{2} )}\frac{1}{\Vert x'-y'\Vert^N}\int_{Q(x',y',x_N)} \frac{|u(x',x_N)-u(y',x_N)|^p}{\Vert x'-y'\Vert^{N-1+sp}}dzdy'dx'dx_N\\
&\leq C\int_{\R^{N-1}}\int_0^{\frac{\eta(x')}{2}}\int_{B(x',\frac{\eta(x')}{2} )}\int_{Q(x',y',x_N)} \frac{|u(x)-u(z)|^p}{\Vert x'-y'\Vert^{2N-1+sp}}dzdy'dx'dx_N\\
&=C\int_\Omega\int_\Omega |u(x)-u(z)|^p \int_{B(x',\frac{\eta(x')}{2} )} \chi_{Q(x',y',x_N)}(z)\frac{1}{\Vert x'
-y'\Vert^{2N-1+sp}} dy' dz dx.
\end{align*}
By (\ref{cubebound}), we have 
$$\frac{4}{3}\Vert x-z\Vert \leq \Vert x'-y'\Vert.$$
Thus, we can simplify to:
\begin{align*}
&\leq C\int_\Omega\int_\Omega |u(x)-u(z)|^p \int_{\R^{N-1}\setminus B(x',\Vert x-z\Vert)} \frac{1}{\Vert x'-y'\Vert^{2N-1+sp}} dy' dz dx\\
&\leq C\int_\Omega\int_\Omega \frac{|u(x)-u(z)|^p }{\Vert x-z\Vert^{N+sp}} dz dx.
\end{align*}
\end{proof}
\begin{theorem}\label{thm:equivdomaingeneral}
Let $\Omega$ be as in (\ref{generalomegawflat}). Let $u\in L^p_{\text{Loc}}(\Omega)$ such that 
\begin{align*}
    \int_{\R^{N-1}}\int_0^{\eta(x')}\int_0^{\eta(x')}\frac{|u(x',x_N)-u(x',y_N)|^p}{|x_N-y_N|^{1+sp}}dy_Ndx_Ndx'&\leq \infty\\
    \int_{\R^{N-1}}\int_0^{\frac{\eta(x')}{2}}\int_{B(x',\frac{\eta(x')}{2} )} \frac{|u(x',x_N)-u(y',x_N)|^p}{\Vert x'-y'\Vert^{N-1+sp}}dy'dx_Ndx'&\leq \infty.
\end{align*}
Then $u\in \dot W^{s,p}(\Omega)$. Furthermore, there is $C=C(N,s,p)>0$ such that
\begin{align*}
    C|u|_{\dot W^{s,p}(\Omega)}^p&\leq \int_{\R^{N-1}}\int_0^{\eta(x')}\int_0^{\eta(x')}\frac{|u(x',x_N)-u(x',y_N)|^p}{|x_N-y_N|^{1+sp}}dy_Ndx_Ndx'\\
    &\tab+\int_{\R^{N-1}}\int_0^{\frac{\eta(x')}{2}}\int_{B(x',\frac{\eta(x')}{2} )} \frac{|u(x',x_N)-u(y',x_N)|^p}{\Vert x'-y'\Vert^{N-1+sp}}dy'dx_Ndx'.
\end{align*}
\end{theorem}
\begin{proof}
We have:
\begin{align*}
\firstline\int_{\R^{N-1}}\int_0^{\eta(x')/8} \int_{B(x',\frac{\eta(x')}{2} )}\int_0^{\eta(y')/8} \frac{|u(x)-u(y)|^p}{\Vert x-y\Vert^{N+sp}} dy dx\\
&=\int_{\R^{N-1}}\int_0^{\eta(x')/8} \int_{B(x',\frac{\eta(x')}{2} )}\int_0^{\eta(y')/8}\frac{|u(x)-u(y)|^p}{\Vert x-y\Vert^{N+sp}} dy_N dy'dx_Ndx'\\
&\leq 2^{p-1} \int_{\R^{N-1}}\int_0^{\eta(x')/8} \int_{B(x',\frac{\eta(x')}{2} )}\int_0^{\eta(y')/8} \frac{|u(x',x_N)-u(y',x_N)|^p}{\Vert x-y\Vert^{N+sp}} dy_N dy'dx_Ndx'\\
&\tab+ 2^{p-1} \int_{\R^{N-1}}\int_0^{\eta(x')/8} \int_{B(x',\frac{\eta(x')}{2} )}\int_0^{\eta(y')/8} \frac{|u(y',x_N)-u(y',y_N)|^p}{\Vert x-y\Vert^{N+sp}} dy_N dy'dx_Ndx'\\
&=:2^{p-1}\mathcal A+2^{p-1}\mathcal B.
\end{align*}
Now, we will look at each part separately. Considering $\mathcal A$:
\begin{align*}
\mathcal A&= \int_{\R^{N-1}}\int_0^{\eta(x')/8} \int_{B(x',\frac{\eta(x')}{2} )}\int_0^{\eta(y')/8} \frac{|u(x',x_N)-u(y',x_N)|^p}{\Vert x-y\Vert^{N+sp}} dy_N dy'dx_Ndx'\\
&=\int_{\R^{N-1}}\int_0^{\eta(x')/8} \int_{B(x',\frac{\eta(x')}{2} )} |u(x',x_N)-u(y',x_N)|^p\int_0^{\eta(y')/8}\frac{1}{\Vert x-y\Vert^{N+sp}} dy_N dy'dx_Ndx'.
\end{align*}
Note that by lemma \ref{lem:slicingfinite}, we get:
$$\leq C\int_{\R^{N-1}}\int_0^{\eta(x')/8} \int_{B(x',\eta(x')/2)} \frac{|u(x',x_N)-u(y',x_N)|^p}{\Vert x'-y'\Vert^{N-1+sp}} dy_N dy'dx_Ndx'.$$
Next, we consider $\mathcal B$:
$$\int_{\R^{N-1}}\int_0^{\eta(x')/8} \int_{B(x',\frac{\eta(x')}{2} )}\int_0^{\eta(y')/8} \frac{|u(y',x_N)-u(y',y_N)|^p}{\Vert x-y\Vert^{N+sp}} dy_N dy'dx_Ndx'.$$
We note that if $\Vert y'-x'\Vert\leq \frac{\eta(x')}{2} $, then we have, using the Lipschitz bound, that $\eta(y')\geq \frac{\eta(x')}{2}$. Therefore, we get:
$$\leq \int_{\mathbb R^{N-1}}\int_0^{\frac{\eta (y')}{4}}\int_0^{\frac{\eta(y')}{4}} |u(y',x_N)-u(y',y_N)|^p \int_{\mathbb R^{N-1}} \frac{1}{\Vert x-y\Vert^{N+sp}}dx'dx_Ndy_Ndy'.$$
By corollary \ref{lem:slicingplane}:
$$\leq \int_{\mathbb R^{N-1}}\int_0^{\frac{\eta (y')}{4}}\int_0^{\frac{\eta(y')}{4}} \frac{|u(y',x_N)-u(y',y_N)|^p}{|x_N-y_N|^{1+sp}} dx_Ndy_Ndy'.$$
\end{proof}

%% file: sections/2.Preliminaries/2.6.MollifierBounds.tex
\subsection{Bounds on the Mollification} \hfill\\
When showing that, for a fixed $g\in L^p_{\text{Loc}}(\mathbb R^{N-1})$, that there is a function $u\in W^{s,p}(\Omega)$ for $\Omega$ as in (\ref{flatomega}) or (\ref{generalomegawflat}) such that $\Tr(u)(\cdot,0)=g(\cdot)$, the natural choice for $u$ is
\begin{equation}\label{eq:mollifiedfunction}
u(x',x_N):=\int_{\mathbb R^{N-1}} \phi(w')g(x'+x_Nw') dw'.
\end{equation}
For $\phi\in C_c^\infty(\mathbb R^{N-1})$ with support in $B(0,1)$ such that $\int_{\mathbb R^{n-1}} \phi dx=1$. Therefore, we naturally wish to show properties of $u$ with respects to $g$.
\begin{lemma}\label{lem:lateralbound}
Let $u$ be as in (\ref{eq:mollifiedfunction}). Then there is $C=C(N,s,p)$ such that:
$$|u(x',x_N)-u(x'+h',x_N)|^p\leq C\frac{\Vert h'\Vert^{p-1}}{x_N^{p+N-1}}\int_0^{\Vert h'\Vert} \int_{B(x'+r\hat h',x_N)} |g(y')-g(x'+r\hat h')|^pdy'dr.$$
\end{lemma}
\begin{proof}
Using the fundamental theorem of calculus and Fubini's theorem, we get
\begin{align*}
u(x'+h',x_N)-u(x',x_N)&=\int_{\R^{N-1}}\frac{1}{x_N^{N-1}} \left(\phi\left(\frac{x'+h'-y'}{x_N}\right)-\phi\left(\frac{x'-y'}{x_N}\right)\right) g(y')dy'dr\\
&=\int_0^{\Vert h'\Vert}\int_{\R^{N-1}}\frac{1}{x_N^{N}} \frac{\partial}{\partial \hat h'}\phi\left(\frac{x'-y'+r\hat h}{x_N}\right) g(y')dy'dr.
\end{align*}
Since the integral of the directional derivative of $\phi$ is $0$, then:
$$=\frac{1}{x_N}\int_0^{\Vert h'\Vert}\int_{\R^{N-1}}\frac{1}{x_N^{N-1}} \frac{\partial}{\partial \hat h'}\phi\left(\frac{x'-y'+r\hat h}{x_N}\right) [g(y')-g(x'+r\hat h')]dy'dr.$$
By H\"older's inequality, we get:
\begin{align*}
|u(x',x_N)-u(x'+h',x_N)|^p&\leq C\frac{\Vert h'\Vert^{p-1}}{x_N^{p+N-1}}\int_0^{\Vert h'\Vert}\int_{\R^{N-1}} \left\vert\frac{\partial}{\partial \hat h'}\phi\left(\frac{x'-y'+r\hat h}{x_N}\right)\right\vert \\
&\tab\times|g(y')-g(x'+r\hat h')|^pdy'dr.
\end{align*}
Since $\phi$ has support in $B(0,1)$, note that the inner term has support for $y'$ in $B(x'+r\hat h,x_N)$, so we rewrite:
$$|u(x',x_N)-u(x'+h',x_N)|^p\leq C\frac{\Vert h\Vert^{p-1}}{x_N^{p+N-1}}\int_0^{\Vert h'\Vert}\int_{B(x'+r\hat h,x_N)} |g(y')-g(x'+r\hat h')|^pdy'dr.$$
\end{proof}
\begin{lemma}\label{verticalbound}
Let $u$ be as in (\ref{eq:mollifiedfunction}). Then there is $C=C(N,s,p)$ such that
$$\int_0^b \int_0^b \frac{|u(x+t)-u(x)|^p}{t^{1+sp}}dtdx\leq C \int_0^{2b} y^{p-sp} |u'(y)|^pdy.$$
\end{lemma}
\begin{proof}
\begin{align*}
\firstline\int_0^b \int_0^b \frac{|u(x+t)-u(x)|^p}{t^{1+sp}}dtdx\\
&=\int_0^b \int_0^b \frac{1}{t^{1+sp}} \left\vert \int_0^{t} u'(x+r)dr\right\vert^pdtdx.
\end{align*}
By Hardy's Inequality:
$$\leq \int_0^b \frac{1}{s^p} \int_0^b (t |u'(x+t)|)^p \frac{dt}{t^{1+sp}}dx.$$
Letting $y=x+t$, we get:
\begin{align*}
&=C\int_0^b \int_x^{x+b} (y-x)^{p-1-sp} |u'(y)|^pdydx\\
&\leq C\int_0^{2b} |u'(y)|^p \int_0^y (y-x)^{p-1-sp}dxdy\\
&=C \int_0^{2b} y^{p-sp} |u'(y)|^pdy.
\end{align*}
\end{proof}

%% file: sections/3.FlatOmega.tex
\section{Trace of Strip Domains}
We now move to the first main result \ref{theorem:flatMT}. We separate this into two parts. We first show properties of the trace, and then show that any functions with these properties are the trace of some homogeneous fractional Sobolev function.   
\subfile{3.FlatOmega/3.1.BoundingTrace}

\subfile{3.FlatOmega/3.2.BoundingMollify}

%% file: sections/3.FlatOmega/3.1.BoundingTrace.tex
\subsection{Bounding Trace} 
\begin{lemma}
For $u\in \dot W^{s,p}(\Omega)$, there is $C=C(N,s,p)>0$ such that:
$$\int_{\R^{N-1}} |u(x',b)-u(x',0)|^p dx'\leq Cb^{sp-1}\vert u\vert_{\dot W^{s,p}(\Omega)}^p.$$
\end{lemma}
\begin{proof}
Let $u\in \dot W^{s,p}(\Omega)$. By theorem \ref{Thm:BoundOnSliceFlat}, we have 
\[\int_{\R^{N-1}} |u(x',\cdot)|_{\dot W^{s,p}((0,b))}^p\;dx'<C|u|_{\dot W^{s,p}(\Omega)}^p.\]
Therefore, for $x'\in \mathbb R^{N-1}$ almost everywhere, $u(x',\cdot)\in \dot W^{s,p}((0,b))$. 
By Morrey's embedding in $W^{s,p}$, for $\alpha=s-1/p$, there is $C=C(s,p)>0$ such that
$$|\bar u(x',\cdot)|_{C^{0,\alpha}((0,b))}\leq C |u(x',\cdot)|_{\dot W^{s,p}((0,b))}.$$
Furthermore, note that since $\bar u$ is H\"older continuous, we can choose values for $\bar u(0)$ and $\bar u(b)$ such that $\bar u:[0,b]\rightarrow \mathbb R$ is still H\"older continuous. Thus, specifically, we have:
$$\frac{|\bar u(x',b)-\bar u(x',0)|}{b^{s-1/p}}\leq C |u(x',\cdot)|_{\dot W^{s,p}([0,b])}.$$
Therefore:
\begin{align*}
\int_{\R^{N-1}} |u(x',b)-u(x',0)|^p dx'&\leq C\int_{\R^{N-1}} |u(x',\cdot)|_{\dot W^{s,p}([0,b])}^p dx'\\
&=C\int_{\R^{N-1}} \int_0^b \int_0^b \frac{|u(x',x_N)-u(x',y_N)|^p}{|x_N-y_N|^{1+sp}} dy_Ndx_Ndx'.
\end{align*}
By theorem \ref{Thm:BoundOnSliceFlat}, we know that this is bounded. 
\end{proof}
Now, the next two theorems bound each side individually. 
\begin{theorem}
For all $u\in \dot W^{s,p}(\Omega)$, for $C=C(N,s,p)>0$, we have that:
$$\int_{\R^{N-1}}\int_{B(x',b)} \frac{|u(x',0)-u(y',0)|^p}{\Vert x'-y'\Vert^{N-2+sp}} dy'dx'\leq C |u|_{\dot W^{s,p}(\Omega)}^p.$$
\end{theorem}
\begin{proof}
Consider
\begin{align*}
\firstline\int_{\R^{N-1}}\int_{B(x',b)} \frac{|u(x',0)-u(y',0)|^p}{\Vert x'-y'\Vert^{N-2+sp}} dy'dx'\\
&=\int_{\R^{N-1}}\int_{B(x',b)}\int_0^{\Vert x'-y'\Vert}\frac{|u(x',0)-u(y',0)|^p}{\Vert x'-y'\Vert^{N-1+sp}} dtdy'dx'\\
&\leq C\int_{\R^{N-1}}\int_{B(x',b)}\int_0^{\Vert x'-y'\Vert}\frac{|u(x',0)-u(x',t)|^p}{\Vert x'-y'\Vert^{N-1+sp}} dtdy'dx'\\
&\tab+C\int_{\R^{N-1}}\int_{B(x',b)}\int_0^{\Vert x'-y'\Vert}\frac{|u(x',t)-u(y',t)|^p}{\Vert x'-y'\Vert^{N-1+sp}} dtdy'dx'\\
&\tab+C\int_{\R^{N-1}}\int_{B(x',b)}\int_0^{\Vert x'-y'\Vert}\frac{|u(y',0)-u(y',t)|^p}{\Vert x'-y'\Vert^{N-1+sp}} dtdy'dx'\\
&:=C\left(\mathcal A_1+\mathcal A_2+\mathcal A_3\right).
\end{align*}
Note that the $\mathcal A_1=\mathcal A_3$ with a change of variables. Therefore, we only need to consider $\mathcal A_1$ and $\mathcal A_2$. 
\begin{align*}
\mathcal A_1&=\int_{\R^{N-1}}\int_{B(x',b)}\int_0^{\Vert x'-y'\Vert}\frac{|u(x',0)-u(x',t)|^p}{\Vert x'-y'\Vert^{N-1+sp}} dtdy'dx'\\
&=\int_{\R^{N-1}}\int_0^b |u(x',0)-u(x',t)|^p \int_{B(x',b)\setminus B(x',t)} \frac{1}{\Vert x'-y'\Vert^{N-1+sp}} dy'dt dx\\
&\leq C\int_{\R^{N-1}}\int_0^b |u(x',0)-u(x',t)|^p \int_{\R^{N-1}\setminus B(x',t)} \frac{1}{\Vert x'-y'\Vert^{N-1+sp}} dy'dt dx\\
&=C\int_{\R^{N-1}} \int_0^b \frac{|u(x',0)-u(x',t)|^p}{t^{sp}} dt dx'.
\end{align*}
Using corollary 1.76 from \cite{fss}, we get:
$$\leq C\int_{\R^{N-1}}\int_0^b \int_0^b \frac{|u(x',x_N)-u(x',y_N)|^p}{|x_N-y_N|^{1+sp}} dy_Ndx_N dx'.$$
This was shown to be bounded in theorem \ref{Thm:BoundOnSliceFlat}. Next, we consider:
$$\int_{\R^{N-1}}\int_{B(x',b)}\int_0^{\Vert x'-y'\Vert}\frac{|u(x',t)-u(y',t)|^p}{\Vert x'-y'\Vert^{N-1+sp}} dtdy'dx'$$
$$\int_{\R^{N-1}}\int_{B(x',b)}\int_0^{b}\frac{|u(x',t)-u(y',t)|^p}{\Vert x'-y'\Vert^{N-1+sp}} dtdy'dx'.$$
We have that this is bounded by theorem \ref{Thm:BoundOnSliceFlat}.
\end{proof}
\begin{theorem}
For all $u\in \dot W^{s,p}(\Omega)$, for $C=C(N,s,p)>0$, we have that:
$$\int_{\R^{N-1}}\int_{\R\setminus B(x',b)} \frac{|u(x',0)-u(y',0)|^p}{\Vert x'-y'\Vert^{N+sp}} dy'dx'\leq \frac{C}{b} |u|_{\dot W^{s,p}(\Omega)}^p.$$
\end{theorem}
\begin{proof}
Consider
\begin{align*}
\firstline\int_{\R^{N-1}}\int_{\R^{N-1}\setminus B(x',b)} \frac{|u(x',0)-u(y',0)|^p}{\Vert x'-y'\Vert^{N+sp}} dy'dx'\\
&=\frac{1}{b}\int_{\R^{N-1}}\int_{\R^{N-1}\setminus B(x',b)}\int_0^b\frac{|u(x',0)-u(y',0)|^p}{\Vert x'-y'\Vert^{N+sp}} dtdy'dx'\\
&\leq \frac{C}{b}\int_{\R^{N-1}}\int_{\R^{N-1}\setminus B(x',b)}\int_0^b\frac{|u(x',0)-u(x',t)|^p}{\Vert x'-y'\Vert^{N+sp}} dtdy'dx'\\
&\tab+\frac{C}{b}\int_{\R^{N-1}}\int_{\R^{N-1}\setminus B(x',b)}\int_0^b\frac{|u(x',t)-u(y',t)|^p}{\Vert x'-y'\Vert^{N+sp}} dtdy'dx'\\
&\tab+\frac{C}{b}\int_{\R^{N-1}}\int_{\R^{N-1}\setminus B(x',b)}\int_0^b\frac{|u(y',0)-u(y',t)|^p}{\Vert x'-y'\Vert^{N+sp}} dtdy'dx'\\
&=:\mathcal A_1+\mathcal A_2+\mathcal A_3.
\end{align*}
Again, we only need to consider $\mathcal A_1$ and $\mathcal A_2$. 
\begin{align*}
\mathcal A_1&=\int_{\R^{N-1}}\int_{\R^{N-1}\setminus B(x',b)}\int_0^b\frac{|u(x',0)-u(x',t)|^p}{\Vert x'-y'\Vert^{N+sp}} dtdy'dx'\\
&=\int_{\R^{N-1}}\int_0^b |u(x',0)-u(x',t)|^p \int_{\R^{N-1}\setminus B(x',b)} \frac{1}{\Vert x'-y'\Vert^{N+sp}} dy'dt dx\\
&\leq C\int_{\R^{N-1}}\int_0^b \frac{|u(x',0)-u(x',t)|^p}{b^{1+sp}} dy'dt dx\\
&\leq \frac{C}{b}\int_{\R^{N-1}} \int_0^b \frac{|u(x',0)-u(x',t)|^p}{t^{sp}} dt dx'.
\end{align*}
Using Corollary 1.76 from \cite{fss}, we get:
$$\leq \frac{C}{b}\int_{\R^{N-1}}\int_0^b \int_0^b \frac{|u(x',x_N)-u(x',y_N)|^p}{|x_N-y_N|^{1+sp}} dy_Ndx_N dx'.$$
This is bounded by theorem \ref{Thm:BoundOnSliceFlat}. Next, we consider:
$$\mathcal A_2=\int_{\R^{N-1}}\int_{\R^{N-1}\setminus B(x',b)}\int_0^b\frac{|u(x',t)-u(y',t)|^p}{\Vert x'-y'\Vert^{N-1+sp}} dtdy'dx'.$$
We have that this is bounded by theorem \ref{Thm:BoundOnSliceFlat}. 
\end{proof}

%% file: sections/3.FlatOmega/3.2.BoundingMollify.tex
\subsection{Bounding the Mollification} \hfill\\
Now, we would like to show that given certain functions on $\mathbb R^{N-1}$, we can find $u\in W^{s,p}(\Omega)$ with trace exactly that function. We first show that given a single function $g$, we can find a function with trace $g$, where we ignore the other side. 
\begin{lemma}\label{lem:onesidelemmaflat}
Let $g\in L^p_{\text{loc}}(\mathbb R^{N-1})$ such that:
$$|g|_{\tilde{W}^{s-1/p,p}_{\leq b}(\mathbb R^{N-1})}<\infty$$
$$|g|_{\tilde{W}^{s+1/p,p}_{\geq b}(\mathbb R^{N-1})}<\infty$$
Then for $u\in L^p_{\text{loc}}(\Omega)$ defined as:
$$g(x',x_N):=\int_{\mathbb R^{N-1}} \phi(w')g(x'+x_Nw') dw',$$
we have that $u\in W^{s,p}(\Omega)$. Furthermore, $g\rightarrow u$ is linear and there is $C=C(N,s,p)>0$ such that
$$|u|_{\dot W^{s,p}(\Omega)}\leq C|g|_{\tilde{W}^{s-1/p,p}_{\leq b}(\mathbb R^{N-1})}+C|g|_{\tilde{W}^{s+1/p,p}_{\geq b}(\mathbb R^{N-1})}.$$
\end{lemma}
\begin{proof}
We will apply theorem \ref{Thm:BoundonSobolevFlat}, so we only consider each term. Primarily, we will consider:
\begin{align*}
\firstline\int_0^b \int_{\R^{N-1}}\int_{B(x',b)} \frac{|u(x',x_N)-u(y',x_N)|^p}{\Vert x'-y'\Vert^{N-1+sp}}dy'dx'dx_N\\
&=\int_0^b \int_{\R^{N-1}}\int_{B(x',x_N)} \frac{|u(x',x_N)-u(y',x_N)|^p}{\Vert x'-y'\Vert^{N-1+sp}}dy'dx'dx_N\\
&\tab+\int_0^b \int_{\R^{N-1}}\int_{B(x',b)\setminus B(x',x_N)} \frac{|u(x',x_N)-u(y',x_N)|^p}{\Vert x'-y'\Vert^{N-1+sp}}dy'dx'dx_N\\
&=\mathcal A_1+\mathcal A_2.
\end{align*}
We will look at these separately. We consider:
\begin{align*}
\mathcal A_{1}=\int_{\mathbb R^{N-1}} \int_0^{b} \int_{B(0,x_N)} \frac{|u(x'+h',x_N)-u(x',x_N)|^p}{\Vert h'\Vert^{N-1+sp}} dh'dx_Ndx'.
\end{align*}
Using lemma \ref{lem:lateralbound}, we get:
\begin{align*}
&\leq \int_{\mathbb R^{N-1}} \int_0^{b} \int_{B(0,x_N)} \frac{|u(x'+h',x_N)-u(x',x_N)|^p}{\Vert h'\Vert^{N-1+sp}} dh'dx_Ndx'\\
&\leq C\int_{\mathbb R^{N-1}} \int_{B(0,b)} \int_{\Vert h'\Vert}^{b} \frac{1}{\Vert h'\Vert^{N-p+sp}x_N^{p+N-1}}\\
&\tab\times\int_0^{\Vert h'\Vert} \int_{B(x'+r\hat h',x_N)} |g(y')-g(x'+r\hat h')|^pdy'drdx_Ndh'dx'.
\end{align*}
Let us use the change of variables $z'=x'+r\hat h'$. We then get:
\begin{align*}
&=C\int_{\mathbb R^{N-1}} \int_{B(0,b)} \int_{\Vert h'\Vert}^{b} \frac{1}{\Vert h'\Vert^{N-p+sp}x_N^{p+N-1}}\int_0^{\Vert h'\Vert} \int_{B(z',x_N)} |g(y')-g(z')|^pdy'drdx_Ndh'dx'\\
&= C\int_{\mathbb R^{N-1}} \int_{B(z',b)} |g(y')-g(z')|^p \int_{\Vert z'-y'\Vert}^{b}\int_{B(0,x_N)} \frac{1}{\Vert h'\Vert^{N-p-1+sp}x_N^{p+N-1}}dh'dx_Ndy'dz'\\
&= C\int_{\mathbb R^{N-1}} \int_{B(z',b)} \frac{|g(y')-g(z')|^p}{\Vert y'-z'\Vert^{N-2+sp}}dy'dz'.
\end{align*}
Next, we consider $\mathcal A_2$. 
\begin{align*}
\mathcal A_2&=\int_{\R^{N-1}} \int_{B(x',b)}\int_0^{\Vert x'-y'\Vert} \frac{|u(x',x_N)-u(y',x_N)|^p}{\Vert x'-y'\Vert^{N-1+sp}}dx_Ndy'dx'\\
&=\int_{\R^{N-1}} \int_{B(x',b)}\int_0^{\Vert x'-y'\Vert} \frac{\left|\int_{\mathbb R^{N-1}}\phi(w')(g(x'+x_Nw')-g(y'+x_Nw'))dw'\right|^p}{\Vert x'-y'\Vert^{N-1+sp}}dx_Ndy'dx'\\
&\leq \int_{\R^{N-1}} \int_{B(x',b)}\int_0^{\Vert x'-y'\Vert} \int_{\mathbb R^{N-1}} \phi(w')\frac{\left|(g(x'+x_Nw')-g(y'+x_Nw')\right|^p}{\Vert x'-y'\Vert^{N-1+sp}}dw'dx_N dy'dx'.
\end{align*}
Using the change of variable $X'=x'+x_Nw'$ and $Y'=y'+x_Nw'$, we get:
\begin{align*}
&=\int_{\R^{N-1}} \int_{B(Y',b)}\int_0^{\Vert X'-Y'\Vert} \int_{\mathbb R^{N-1}} \phi(w')\frac{\left|g(X')-g(Y')\right|^p}{\Vert X'-Y'\Vert^{N-1+sp}}dw'dx_N dY'dX'\\
&=\int_{\R^{N-1}} \int_{B(Y',b)}\frac{\left|(g(X')-g(Y')\right|^p}{\Vert X'-Y'\Vert^{N-2+sp}}dY'dX'
\end{align*}
Therefore, we have achieved the local bound. Now, we will look at the far bound. 
\begin{align*}
\firstline \int_0^b \int_{\R^{N-1}} \int_{\R^{N-1}\setminus B(x',b)} \frac{|u(x',x_N)-u(y',x_N)|^p}{\Vert x'-y'\Vert^{N+sp}}dy'dx'dx_N\\
&=\int_0^b\int_{\R^{N-1}} \int_{\R^{N-1}\setminus B(x',b)} \frac{\left|\int_{\mathbb R^{N-1}}\phi(w')(g(x'+x_Nw')-g(y'+x_Nw')dw'\right|^p}{\Vert x'-y'\Vert^{N+sp}}dy'dx'dx_N\\
&\leq \int_0^b\int_{\R^{N-1}} \int_{\R^{N-1}\setminus B(x',b)} \int_{\mathbb R^{N-1}} \phi(w')\frac{\left|(g(x'+x_Nw')-g(y'+x_Nw')\right|^p}{\Vert x'-y'\Vert^{N+sp}}dw'dy'dx'dx_N.
\end{align*}
Using the change of variable $X'=x'+x_Nw'$ and $Y'=y'+x_Nw'$, we get:
\begin{align*}
&=\int_0^b\int_{\R^{N-1}} \int_{\R^{N-1}\setminus B(x',b)} \int_{\mathbb R^{N-1}} \phi(w')\frac{\left|(g(X')-g(Y')\right|^p}{\Vert X'-Y'\Vert^{N+sp}}dw'dY'dX'dx_N\\
&=b\int_{\R^{N-1}} \int_{\R^{N-1}\setminus B(x',b)}\frac{\left|(g(X')-g(Y')\right|^p}{\Vert X'-Y'\Vert^{N+sp}}dY'dX'.
\end{align*}
Now, lastly, we need a bound along the $(0,b)$ direction.
\begin{align*}
\firstline\int_{\mathbb R^{N-1}} \int_0^{b} \int_0^{b} \frac{|u(x',x_N)-u(x',y_N)|^p}{|x_N-y_N|^{1+sp}}dy_Ndx_Ndx'\\
&=\int_{\mathbb R^{N-1}} \int_0^{b} \int_0^{b} \frac{|u(x',x_N+t)-u(x',x_N)|^p}{t^{1+sp}}dtdx_Ndx'.
\end{align*}
By lemma \ref{verticalbound}, we have:
$$\leq C\int_{\mathbb R^{N-1}} \int_0^{2b} x_N^{p-sp} |\partial_N u(x',x_N)|^pdx_Ndx'.$$
Now, we use the inequality
$$\vert \partial_N u(x)\vert \leq C \frac{1}{x_N^N} \int_{B(x',x_N)}|g(y')-g(x')|dy'$$
to get
$$\leq C \int_{\mathbb R^{N-1}} \int_0^{2b}  x_N^{p-sp}\left(\int_{B(x',x_N)} \frac{1}{x_N^{N}} |g(y')-g(x')|dy'\right)^pdx_Ndx'.$$
By H\"older's inequality:
\begin{align*}
&\leq C \int_{\mathbb R^{N-1}} \int_0^{2b}  x_N^{p-sp-Np+(N-1)(p-1)}\int_{B(x',x_N)} |g(y')-g(x')|^pdy'dx_Ndx'\\
&=C\int_{\mathbb R^{N-1}} \int^{2b} |g(y')-g(x')|^p \int_{\Vert x'-y'\Vert}^{2b} \frac{1}{x_N^{N-1+sp}}dx_Ndy'dx'\\
&\leq C\int_{\mathbb R^{N-1}}\int_{B(x',2b)} \frac{|g(x')-g(y')|^p}{\Vert x'-y'\Vert^{sp}}dy'dx'.
\end{align*}
\end{proof}
Now, lastly, we consider the general case in theorem~\ref{theorem:flatMT}.  
\begin{proof}
First, note that if we have $g_+$ and $g_-$ satisfying the condition, let $u_+$ and $u_-$ be their corresponding functions as generated by lemma~\ref{lem:onesidelemmaflat}. Consider $h(x'):=g_+(x')-u_-(x',b)$. Note that:
\begin{align*}
\int_{\R^{N-1}} |h(x')|^pdx &\leq C\int_{\R^{N-1}} |g_+(x')-g_-(x')|^p dx'+C\int_{\R^{N-1}} |u_-(x',b)-g_-(x')|^p dx'\\
&=C\int_{\R^{N-1}} |g_+(x')-g_-(x')|^p dx'+C\int_{\R^{N-1}} |u_-(x',b)-u_-(x',0)|^p dx'.
\end{align*}
We are given that the first term is finite, and due to $u\in \dot W^{s,p}(\Omega)$, by lemma \ref{Thm:BoundOnSliceFlat}, we have the second term is finite. Therefore, if we can find a function $v\in \dot W^{s,p}(\Omega)$ such that $v(x',b)=h(x')$, and $v(x',0)=0$, then we would have that $f(x)=v+u_-$ would satisfy our requirements. Let $\psi\in C^\infty([0,b])$ such that $\psi(b)=0$ and $\psi(0)=1$. Let $v(x)=\psi(x_N)u_h(x',b-x_N)$. Note that other than in the $(0,b)$ direction, the bounds from the previous part hold. Therefore, we need to consider:
\begin{align*}
\firstline\int_{\R^{N-1}}\int_0^b \int_0^b \frac{|\Delta_{t,N} v(x)|^p}{t^{1+sp}} dx_Ndtdx'\\
&=\int_{\R^{N-1}}\int_0^b \int_0^b \frac{|\psi(x_N+t)u_h(x',x_N+t)-\psi(x_N)u_h(x',x_N)|^p}{t^{1+sp}} dx_Ndtdx'\\
&\leq C\int_{\R^{N-1}}\int_0^b \int_0^b \frac{|\psi(x_N+t)u_h(x',x_N+t)-\psi(x_N+t)u_h(x',x_N)|^p}{t^{1+sp}} dx_Ndtdx'\\
&\tab+C\int_{\R^{N-1}}\int_0^b \int_0^b \frac{|\psi(x_N+t)u_h(x',x_N)-\psi(x_N)u_h(x',x_N)|^p}{t^{1+sp}} dx_Ndtdx'\\
&=C\mathcal B_1+C\mathcal B_2.
\end{align*}
We would have:
\begin{align*}
\mathcal B_2&=\int_{\R^{N-1}}\int_0^b \int_0^b |\psi(x_N+t)|\frac{|u_h(x',x_N+t)-u_h(x',x_N)|^p}{t^{1+sp}} dx_Ndtdx'\\
&\leq \int_{\R^{N-1}}\int_0^b \int_0^b \frac{|u_h(x',x_N+t)-u_h(x',x_N)|^p}{t^{1+sp}} dx_Ndtdx'.
\end{align*}
We have already shown this is bounded. Next, consider:
\begin{align*}
\mathcal B_1=\int_{\R^{N-1}}\int_0^b \int_0^b \frac{|\psi(x_N+t)-\psi(x_N)|^p}{t^{1+sp}} |u_h(x)|^p dx_Ndtdx'.
\end{align*}
Since $\psi$ has bounded derivative, we get
\begin{align*}
\leq C\int_{\R^{N-1}}\int_0^b \int_0^b \frac{t^p}{t^{1+sp}} |u_h(x)|^p dx_Ndtdx'
\leq C\int_{\R^{N-1}}\int_0^b |u_h(x)|^p dx_Ndtdx'.
\end{align*}
Since $\phi\in W^{s,p}([0,b])$ and $h\in L^p$, we have this is finite. 
\end{proof}

%% file: sections/4.GeneralOmega.tex
\section{Trace of Strip-like Domains}
\subfile{4.GeneralOmega/4.1.BoundingTrace}

\subfile{4.GeneralOmega/4.2.BoundingMollify}

%% file: sections/4.GeneralOmega/4.1.BoundingTrace.tex
\subsection{Bounding Trace} \hfill\\
Now we will apply the previous part to bound the trace. First, we will bound the difference of the two sides. 
\begin{theorem}[Bounding the Difference]
For $u\in \dot W^{s,p}(\Omega)$, there is $C=C(N,s,p)>0$ such that:
$$\int_{\R^{N-1}} \frac{|u(x',b)-u(x',0)|^p}{
\eta(x')^{sp-1}}dx'\leq C\vert u\vert_{\dot W^{s,p}(\Omega)}^p$$
\end{theorem}
\begin{proof}
Let $u\in \dot W^{s,p}(\Omega)$. By theorem $4.1$, we have:
$$\int_{\R^{N-1}}\int_0^{\eta(x')}\int_0^{\eta(x')}\frac{|u(x',x_N)-u(x',y_N)|^p}{|x_N-y_N|^{1+sp}}dy_Ndx_Ndx'\leq C\vert u\vert_{\dot W^{s,p}(\Omega)}^p$$
Therefore, for $x'\in \mathbb R^{N-1}$ almost everywhere, $u(x',\cdot)\in W^{s,p}((0,\eta(x'))$. 
By Morrey's embedding in $W^{s,p}$, for $\alpha=s-1/p$ there is $C=C(s,p)>0$ such that
$$|\bar u(x',\cdot)|_{C^{0,\alpha}((0,\eta(x')))}\leq C |u(x',\cdot)|_{\dot W^{s,p}((0,\eta(x')))}$$
Furthermore, note that since $\bar u$ is H\"older continuous, we can choose values for $\bar u(0)$ and $\bar u(\eta(x'))$ such that $\bar u:[0,\eta(x')]\rightarrow \mathbb R$ is still H\"older continuous, Thus, specifically, we have:
$$\frac{|\bar u(x',b)-\bar u(x',0)|}{\eta(x')^{s-1/p}}\leq C |u(x',\cdot)|_{\dot W^{s,p}([0,\eta(x')])}$$
Therefore:
\begin{align*}
\int_{\R^{N-1}} \frac{|\bar u(x',b)-\bar u(x',0)|^p}{\eta(x')^{sp-1}} dx'&\leq C\int_{\R^{N-1}} |u(x',\cdot)|_{\dot W^{s,p}([0,\eta(x')])}^p dx'\\
&=C\int_{\R^{N-1}} \int_0^{\eta(x')} \int_0^{\eta(x')} \frac{|u(x',x_N)-u(x',y_N)|^p}{|x_N-y_N|^{1+sp}} dy_Ndx_Ndx'\\
\end{align*}
By Theorem \ref{thm:equivdomaingeneral}, we know that this is bounded.
\end{proof}
Now, the next two theorems bound each side individually. 
\begin{theorem}
For $u\in \dot W^{s,p}(\Omega)$, for $C=C(N,s,p )>0$, we have that:
$$|u(\cdot,0)|_{\tilde{W}^{s-1/p,p}_{\leq \eta/2}(\mathbb R^{N-1})}\leq C |u|_{\dot W^{s,p}(\Omega)}$$
\end{theorem}
\begin{proof}
Consider:
\begin{align*}
\firstline\int_{\R^{N-1}}\int_{B(x',\frac{\eta(x')}{2} )} \frac{|u(x',0)-u(y',0)|^p}{\Vert x'-y'\Vert^{N-2+sp}} dy'dx'\\
&=\int_{\R^{N-1}}\int_{B(x',\frac{\eta(x')}{2} )}\int_0^{\Vert x'-y'\Vert}\frac{|u(x',0)-u(y',0)|^p}{\Vert x'-y'\Vert^{N-1+sp}} dtdy'dx'\\
&\leq C\int_{\R^{N-1}}\int_{B(x',\frac{\eta(x')}{2} )}\int_0^{\Vert x'-y'\Vert}\frac{|u(x',0)-u(x',t)|^p}{\Vert x'-y'\Vert^{N-1+sp}} dtdy'dx'\\
&\tab+C\int_{\R^{N-1}}\int_{B(x',\frac{\eta(x')}{2} )}\int_0^{\Vert x'-y'\Vert}\frac{|u(x',t)-u(y',t)|^p}{\Vert x'-y'\Vert^{N-1+sp}} dtdy'dx'\\
&\tab+C\int_{\R^{N-1}}\int_{B(x',\frac{\eta(x')}{2} )}\int_0^{\Vert x'-y'\Vert}\frac{|u(y',0)-u(y',t)|^p}{\Vert x'-y'\Vert^{N-1+sp}} dtdy'dx'\\
&=C(\mathcal A_1+\mathcal A_2+\mathcal A_3)
\end{align*}
First:
\begin{align}
\mathcal A_1&=\int_{\R^{N-1}}\int_{B(x',\frac{\eta(x')}{2} )}\int_0^{\Vert x'-y'\Vert}\frac{|u(x',0)-u(x',t)|^p}{\Vert x'-y'\Vert^{N-1+sp}} dtdy'dx'\label{A1}\\
&=\int_{\R^{N-1}}\int_0^{\eta(x')} |u(x',0)-u(x',t)|^p \int_{B(x',\frac{\eta(x')}{2} )\setminus B(x',t)} \frac{1}{\Vert x'-y'\Vert^{N-1+sp}} dy'dt dx\nonumber\\
&\leq C\int_{\R^{N-1}}\int_0^{\eta(x')} |u(x',0)-u(x',t)|^p \int_{\R^{N-1}\setminus B(x',t)} \frac{1}{\Vert x'-y'\Vert^{N-1+sp}} dy'dt dx\nonumber\\
&=C\int_{\R^{N-1}} \int_0^{\eta(x')}  \frac{|u(x',0)-u(x',t)|^p}{t^{sp}} dt dx'\nonumber
\end{align}
Using Corollary 1.76 from \cite{fss}, we get:
$$\leq C\int_{\R^{N-1}}\int_0^{\eta(x')} \int_0^{\eta(x')} \frac{|u(x',x_N)-u(x',y_N)|^p}{|x_N-y_N|^{1+sp}} dy_Ndx_N dx'$$
This is bounded by Theorem $4.1$. Next, we consider:
\begin{align*}
\mathcal A_2&=\int_{\R^{N-1}}\int_{B(x',\frac{\eta(x')}{2} )}\int_0^{\Vert x'-y'\Vert}\frac{|u(x',t)-u(y',t)|^p}{\Vert x'-y'\Vert^{N-1+sp}} dtdy'dx'\\
&\leq \int_{\R^{N-1}}\int_{B(x',\frac{\eta(x')}{2} )}\int_0^{\frac{\eta(x')}{2} }\frac{|u(x',t)-u(y',t)|^p}{\Vert x'-y'\Vert^{N-1+sp}} dtdy'dx'
\end{align*}
We have that this is bounded by Theorem $4.2$. Lastly, we have:
\begin{align*}
\mathcal A_3&=\int_{\R^{N-1}}\int_{B(x',\frac{\eta(x')}{2} )}\int_0^{\Vert x'-y'\Vert}\frac{|u(y',0)-u(y',t)|^p}{\Vert x'-y'\Vert^{N-1+sp}} dtdy'dx'\\
&=\int_{\R^{N-1}}\int_{\R^{N-1}}\chi_{B(x',\frac{\eta(x')}{2} )}(y')\int_0^{\Vert x'-y'\Vert}\frac{|u(y',0)-u(y',t)|^p}{\Vert x'-y'\Vert^{N-1+sp}} dtdy'dx'
\end{align*}
Note specifically that if $y'\in B(x',\frac{\eta(x')}{2} )$, then we have that $\Vert y'-x'\Vert\leq \frac{\eta(x')}{2} $. Therefore, using the Lipschitz bound, $\eta(y')\geq \eta(x')-\Vert x'-y'\Vert\geq \frac{\eta(x')}{2}$. Therefore, we have $\Vert y'-x'\Vert\leq \eta(y') $. Thus, get the bound: 
\begin{align*}
\leq \int_{\R^{N-1}}\int_{B(y',\eta(y'))}\int_0^{\Vert x'-y'\Vert}\frac{|u(y',0)-u(y',t)|^p}{\Vert x'-y'\Vert^{N-1+sp}} dtdx'dy'
\end{align*}
We apply the work from (\ref{A1})
\end{proof}
\begin{theorem}
For $u\in \dot W^{s,p}(\Omega)$, for $C=C(N,s,p)>0$, we have that:
$$|u(\cdot,0)|_{\tilde{W}^{s+1/p,p}_{\geq \eta/2}(\mathbb R^{N-1})}\leq C |u|_{\dot W^{s,p}(\Omega)}$$
\end{theorem}
\begin{proof}
We have:
\begin{align*}
\firstline\int_{\mathbb R^{N-1}}\int_{\mathbb R^{N-1}\setminus B(x',\frac{\eta(x')}{2} )} \eta(x')\eta(y')\frac{|u(x',0)-u(y',0)|^p}{\Vert x'-y'\Vert^{N+sp}} dy'dx'\\
&=\int_{\mathbb R^{N-1}}\int_{\mathbb R^{N-1}\setminus B(x',\frac{\eta(x')}{2} )} \int_0^{\eta(x')}\int_0^{\eta(y')}\frac{|u(x',0)-u(y',0)|^p}{\Vert x'-y'\Vert^{N+sp}} dy_Ndx_Ndy'dx'\\
&=C\int_{\mathbb R^{N-1}}\int_{\mathbb R^{N-1}\setminus B(x',\frac{\eta(x')}{2} )} \int_0^{\eta(x')}\int_0^{\eta(y')}\frac{|u(x',x_N)-u(x',0)|^p}{\Vert x'-y'\Vert^{N+sp}} dy_Ndx_Ndy'dx'\\
&\tab+C\int_{\mathbb R^{N-1}}\int_{\mathbb R^{N-1}\setminus B(x',\frac{\eta(x')}{2} )} \int_0^{\eta(x')}\int_0^{\eta(y')}\frac{|u(x',x_N)-u(y',y_N)|^p}{\Vert x'-y'\Vert^{N+sp}} dy_Ndx_Ndy'dx'\\
&\tab+C\int_{\mathbb R^{N-1}}\int_{\mathbb R^{N-1}\setminus B(x',\frac{\eta(x')}{2} )} \int_0^{\eta(x')}\int_0^{\eta(y')}\frac{|u(y',0)-u(y',y_N)|^p}{\Vert x'-y'\Vert^{N+sp}} dy_Ndx_Ndy'dx'\\
&=C(\mathcal A_1+\mathcal A_2+\mathcal A_3)
\end{align*}
We bound each of these individually. First, we have:
\begin{align}
\mathcal A_1&=\int_{\mathbb R^{N-1}}\int_{\mathbb R^{N-1}\setminus B(x',\frac{\eta(x')}{2} )} \int_0^{\eta(x')}\int_0^{\eta(y')}\frac{|u(x',x_N)-u(x',0)|^p}{\Vert x'-y'\Vert^{N+sp}} dy_Ndx_Ndy'dx'\label{A12}\\
&=\int_{\mathbb R^{N-1}}\int_0^{\eta(x')}|u(x',x_N)-u(x',0)|^p\int_{\mathbb R^{N-1}\setminus B(x',\frac{\eta(x')}{2} )}\frac{\eta(y')}{\Vert x'-y'\Vert^{N+sp}} dy'dx_Ndx'\nonumber\\
&\leq \int_{\mathbb R^{N-1}}\int_0^{\eta(x')}|u(x',x_N)-u(x',0)|^p\int_{\mathbb R^{N-1}\setminus B(x',\frac{\eta(x')}{2} )}\frac{\eta(x')+ \Vert x'-y'\Vert}{\Vert x'-y'\Vert^{N+sp}} dy'dx_Ndx'\nonumber\\
&\leq \int_{\mathbb R^{N-1}}\int_0^{\eta(x')}|u(x',x_N)-u(x',0)|^p\int_{\mathbb R^{N-1}\setminus B(x',\frac{\eta(x')}{2} )}\frac{2 \Vert x'-y'\Vert+ \Vert x'-y'\Vert}{\Vert x'-y'\Vert^{N+sp}} dy'dx_Ndx'\nonumber\\
&=C\int_{\mathbb R^{N-1}}\int_0^{\eta(x')}|u(x',x_N)-u(x',0)|^p\int_{\mathbb R^{N-1}\setminus B(x',\frac{\eta(x')}{2} )}\frac{1}{\Vert x'-y'\Vert^{N-1+sp}} dy'dx_Ndx'\nonumber\\
&\leq C\int_{\mathbb R^{N-1}}\int_0^{\eta(x')}\frac{|u(x',x_N)-u(x',0)|^p}{\eta(x')^{sp}} dx_Ndx'\nonumber\\
&\leq C\int_{\mathbb R^{N-1}}\int_0^{\eta(x')}\frac{|u(x',x_N)-u(x',0)|^p}{x_N^{sp}} dx_Ndx'\nonumber
\end{align}

Using Corollary $1.76$ from \cite{fss}, we get:
$$\leq C\int_{\mathbb R^{N-1}} \int_0^{\eta(x')}\int_0^{\eta(x')} \frac{|u(x',x_N)-u(x',y_N)|^p}{|x_N-y_N|^{1+sp}} dy_Ndx_Ndx'$$
We can notice that with the same proof used to derive $5.2.1$:
$$\mathcal A_3\leq \int_{\mathbb R^{N-1}}\int_{\mathbb R^{N-1}\setminus B(y',\eta(y') )} \int_0^{\eta(x')}\int_0^{\eta(y')}\frac{|u(x',x_N)-u(x',0)|^p}{\Vert x'-y'\Vert^{N+sp}} dy_Ndx_Ndy'dx'$$
Now, using similar reasoning to (\ref{A12}), we get a bound. Lastly, we note that with $\mathcal A_2$, we get:$$\mathcal A_2=\int_{\mathbb R^{N-1}}\int_{\mathbb R^{N-1}\setminus B(x',\frac{\eta(x')}{2} )} \int_0^{\eta(x')}\int_0^{\eta(y')}\frac{|u(x',x_N)-u(y',y_N)|^p}{\Vert x'-y'\Vert^{N+sp}} dy_Ndx_Ndy'dx'$$
We want to compare $\Vert x'-y'\Vert$ with $\Vert x-y\Vert$. We note that $\Vert x-y\Vert \leq \Vert x'-y'\Vert+|x_N-y_N|$. We then note that $x_N-y_N$ is bounded by $\max(\eta(x'),\eta(y'))$.  Lastly, we note $\eta(y')$ is bounded by $\eta(x')+ \Vert x'-y'\Vert$. Since we have $\frac{\eta(x')}{2} \leq \Vert x'-y'\Vert$, we get an overall bound of $3 \Vert x'-y'\Vert$. Therefore, we get:
$$\leq C\int_{\mathbb R^{N-1}}\int_{\mathbb R^{N-1}\setminus B(x',\frac{\eta(x')}{2} )} \int_0^{\eta(x')}\int_0^{\eta(y')}\frac{|u(x',x_N)-u(y',y_N)|^p}{\Vert x-y\Vert^{N+sp}} dy_Ndx_Ndy'dx'$$
$$\leq |u|_{\dot W^{s,p}(\Omega)}^p$$
\end{proof}

%% file: sections/4.GeneralOmega/4.2.BoundingMollify.tex
\subsection{Bounding the Mollification} \hfill\\
Now, we would like to show that given certain functions on $\mathbb R^{N-1}$, we can find $u\in W^{s,p}(\Omega)$ with trace exactly that function. The end result we would like to show is the following theorem. 

We will only consider one side for now. This gives us we have the correct condition on each side, and we need only verify that the cross-side condition is correct. 
\begin{lemma}\label{lem:onesidegeneral}
Let $\Omega$ be as in (\ref{generalomegawflat}). Let $g\in L^p_{\text{loc}}(\mathbb R^{N-1})$ such that:
\begin{align*}
|g|_{\tilde{W}^{s-1/p,p}_{\leq \eta/2}(\R^{N-1})}&<\infty\\
|g|_{\tilde{W}^{s+1/p,p}_{\geq \eta/2}(\R^{N-1})}&<\infty
\end{align*}
Then for $u\in L^p_{\text{loc}}(\Omega)$ defined as:
$$u(x',x_N):=\int_{\mathbb R^{N-1}} \phi(w')g(x'+\frac{x_N}{8}w') dw'$$
We have that $u\in W^{s,p}(\Omega)$. Furthermore, $g\rightarrow u$ is linear and there is $C=C(N,s,p,L)>0$ such that
\begin{align*}
C|u|_{\dot W^{s,p}(\Omega')}&\leq |g|_{\tilde{W}^{s-1/p,p}_{\leq \eta/2}(\R^{N-1})}+|g|_{\tilde{W}^{s+1/p,p}_{\geq \eta/2}(\R^{N-1})}
\end{align*}
\end{lemma}
\begin{proof}
We will show, for:
$$u(x',x_N):=\int_{\mathbb R^{N-1}} \phi(w')g(x'+x_Nw') dw',$$
$u\in W^{s,p}(\Omega')$ for $\Omega':=\{(x',x_N):0<x_N<\eta(x')/8\}$; however, using lemma \ref{lem:dilation}, we can turn this function into one over $\Omega$. 
\par We consider:
\begin{align*}
\int_{\Omega'}\int_{\Omega'} \frac{|u(x)-u(y)|^p}{\Vert x-y\Vert^{N+sp}}dydx
&=\int_{\Omega'}\int_{B(x',\eta(x')/2) }\int_0^{\eta(y')/8} \frac{|u(x)-u(y)|^p}{\Vert x-y\Vert^{N+sp}}dy_Ndy'dx\\
&\tab+\int_{\Omega'}\int_{\mathbb R^{N-1}\setminus B(x',\eta(x')/2)} \int_0^{\eta(y')/8} \frac{|u(x)-u(y)|^p}{\Vert x-y\Vert^{N+sp}}dy_Ndy'dx\\
&=:\mathcal A+\mathcal B.
\end{align*}
We bound each term separately. First, consider $\mathcal B$:
\begin{align*}
\mathcal B&=\int_{\Omega'}\int_{\mathbb R^{N-1}\setminus B(x',\frac{\eta(x')}{2} )} \int_0^{\eta(y')/8} \frac{|u(x)-u(y)|^p}{\Vert x-y\Vert^{N+sp}}dy_Ndy'dx\\
&=\int_{\Omega'}\int_{\mathbb R^{N-1}\setminus B(x',\frac{\eta(x')}{2})} \int_0^{\eta(y')/8} \frac{\left\vert\int_{\mathbb R^{N-1}} \phi(z')(g(x'+x_Nz')-g(y'+y_Nz'))dz'\right\vert^p}{\Vert x-y\Vert^{N+sp}}dy_Ndy'dx.
\end{align*}
By H\"older's inequality:
\begin{align*}
\leq C\int_{\Omega'}\int_{\mathbb R^{N-1}\setminus B(x',\frac{\eta(x')}{2})} \int_0^{\eta(y')/8} \frac{\int_{\mathbb R^{N-1}} \phi(z')|g(x'+x_Nz')-g(y'+y_Nz')|^pdz'}{\Vert x-y\Vert^{N+sp}}dydx\\
\leq C\int_{\Omega'}\int_{\mathbb R^{N-1}\setminus B(x',\frac{\eta(x')}{2} )} \int_0^{\eta(y')/8} \int_{B(0,1)}\frac{ |g(x'+x_Nz')-g(y'+y_Nz')|^p}{\Vert x-y\Vert^{N+sp}}dz'dydx.
\end{align*}
We note that since $\eta$ is $1$-Lipschitz continuous
$$\Vert x'+x_Nz'-y'-y_Nz'\Vert\leq \Vert x'-y'\Vert +\Vert  x_Nz'\Vert +\Vert y_Nz'\Vert\leq \Vert x'-y'\Vert +\frac{\eta(x')}{8}+\frac{\eta(y')}{8}$$
$$\leq \Vert x'-y'\Vert + \frac{\eta(x')}{8}+\frac{\eta(x')+ \Vert x'-y'\Vert}{8}=\frac{9}{8}\Vert x'-y'\Vert+\frac{1}{4}\eta(x')\leq \frac{9}{8}\Vert x'-y'\Vert + \frac{1}{2}\Vert x'-y'\Vert=C\Vert x'-y'\Vert.$$
Therefore, we get:
\begin{align*}
&\leq C\int_{\Omega'}\int_{\mathbb R^{N-1}\setminus B(x',\frac{\eta(x')}{2} )} \int_0^{\eta(y')/8} \int_{B(0,1)}\frac{ |g(x'+x_Nz')-g(y'+y_Nz')|^p}{\Vert x'+x_Nz'-y'-y_Nz'\Vert^{N+sp}}dz'dy_Ndy'dx\\
&= C\int_{B(0,1)}\int_0^\infty \int_{\mathbb R^{N-1}} \int_0^\infty \int_{\mathbb R^{N-1}\setminus B(x',\frac{\eta(x')}{2} )} \chi_{(0,\eta(x')/8)}(x_N)\chi_{(0,\eta(y')/8)}(y_N)\\
&\tab\times \frac{ |g(x'+x_Nz')-g(y'+y_Nz')|^p}{\Vert x'+x_Nz'-y'-y_Nz'\Vert^{N+sp}}dy_Ndy'dxdz'.
\end{align*}
Using a change of variables $X'=x'+x_Nz'$ and $Y'=y'+y_Nz'$
\begin{align*}
&= C\int_{B(0,1)}\int_0^\infty \int_{\mathbb R^{N-1}} \int_0^\infty \int_{\mathbb R^{N-1}\setminus B(X'-x_Nz',\frac{\eta(X'-x_Nz')}{2} )} \chi_{\left(0,\frac{\eta(X'- x_Nz')}{8}\right)}(x_N)\chi_{\left(0,\frac{\eta(Y'-y_Nz')}{8}\right)}(y_N)\\
&\tab\times\frac{ |g(X')-g(Y')|^p}{\Vert X'-Y'\Vert^{N+sp}}dY'dy_NdX'dx_Ndz'.
\end{align*}
\normalsize
We first note that $\eta(X'-x_Nz')\leq \eta(X')+x_N $ since $\eta$ is $1$-Lipschitz continuous. Using the indicator function, we know $x_N\leq \eta(X'-x_Nz')/8$; therefore, we have:
$$\eta(X'-x_Nz')\leq \eta(X')+\eta(X'-x_Nz')/8.$$
We rewrite to
$$\eta(X'-x_Nz')\leq \frac{8}{7}\eta(X').$$
We can do the same for $Y'$. Similarly, we get a lower bound for it as:
$$\eta(X'-x_Nz')\geq \eta(X')-\eta(X'-x_Nz')/8.$$
That is,
$$\eta(X'-x_Nz')\geq \frac{8}{9}\eta(X').$$
Furthermore, we note that:
$$\Vert x_Nz'\Vert \leq \frac{1}{8}\eta(X'-x_Nz')\leq \frac{1}{7}\eta(X').$$
Therefore, we can consider $w'\in B\left(X',\left(\frac{4}{9}-\frac{1}{7}\right)\eta(X')\right)$; then:
$$\Vert w'-X'+x_Nz'\Vert\leq \Vert w'-X'\Vert+\Vert x_Nz'\Vert\leq \left(\frac{4}{9}-\frac{1}{7}\right)\eta(X')+\frac{1}{7}\eta(X')=\frac{4}{9}\eta(X')\leq \frac{\eta(X'-x_Nz')}{2}.$$
Therefore
$$B\left(X',\left(\frac{4}{9}-\frac{1}{7}\right)\eta(X')\right)\subseteq B\left(X'-x_Nz',\frac{\eta(X'-x_Nz')}{2} \right).$$
The constant comes out to $19/63$. Thus, we can simplify the above expression to:
\begin{align*}
&\leq C\int_{B(0,1)}\int_0^\infty \int_{\mathbb R^{N-1}} \int_0^\infty \int_{\mathbb R^{N-1}\setminus B(X',19\eta(X')/63)} \chi_{(0,\eta(X')/7)}(x_N)\chi_{(0,\eta(Y')/7)}(y_N)\\ 
&\tab\times\frac{ |g(X')-g(Y')|^p}{\Vert X'-Y'\Vert^{N+sp}}dY'dy_NdX'dx_Ndz'\\
&=C\int_{\mathbb R^{N-1}}  \int_{\mathbb R^{N-1}\setminus B(X',19\eta(X')/63)} \eta(X')\eta(Y') \frac{ |g(X')-g(Y')|^p}{\Vert X'-Y'\Vert^{N+sp}}dY'dX'.
\end{align*}
This is bounded by lemma \ref{lem:decreaser}. Therefore, we have bounded the $\mathcal B$ term. Consider the $\mathcal A$ term; we have:
$$=\int_{\Omega'}\int_{B(x',\eta(x')/2) }\int_0^{\eta(y')/8} \frac{|u(x)-u(y)|^p}{\Vert x-y\Vert^{N+sp}}dydx.$$
By our lemma \ref{thm:equivdomaingeneral}, we have a bound of:
\begin{align*}
&\leq C\int_{\mathbb R^{N-1}} \int_0^{\eta(x')/8} \int_{B(x',\frac{\eta(x')}{2})} \frac{|u(x',x_N)-u(y',x_N)|^p}{\Vert x'-y'\Vert^{N-1+sp}}dy'dx_Ndx'\\
&\tab+ C\int_{\mathbb R^{N-1}} \int_0^{\eta(x')/4} \int_0^{\eta(x')/4} \frac{|u(x',x_N)-u(y',y_N)|^p}{|x_N-y_N|^{1+sp}}dy_Ndx_Ndx'\\
&=:C(\mathcal A_1+\mathcal A_2).
\end{align*}
We first look at:
\begin{align*}
\mathcal A_1&:=\int_{\mathbb R^{N-1}} \int_0^{\eta(x')/8} \int_{B(0,\eta(x')/2)} \frac{|u(x'+h',x_N)-u(x',x_N)|^p}{\Vert h'\Vert^{N-1+sp}} dh'dx_Ndx'\\
&=\int_{\mathbb R^{N-1}} \int_0^{\eta(x')/8} \int_{B(0,x_N)} \frac{|u(x'+h',x_N)-u(x',x_N)|^p}{\Vert h'\Vert^{N-1+sp}} dh'dx_Ndx'\\
&\tab+\int_{\mathbb R^{N-1}} \int_0^{\eta(x')/8} \int_{B(0,\eta(x')/2)\setminus B(0,x_N)} \frac{|u(x'+h',x_N)-u(x',x_N)|^p}{\Vert h'\Vert^{N-1+sp}} dh'dx_Ndx'\\
&=:\mathcal A_{1,1}+\mathcal A_{1,2}.
\end{align*}
We then consider:
\begin{align*}
\mathcal A_{1,2}&=\int_{\mathbb R^{N-1}} \int_0^{\eta(x')/8} \int_{B(0,\eta(x')/2)\setminus B(0,x_N)} \frac{|u(x'+h',x_N)-u(x',x_N)|^p}{\Vert h'\Vert^{N-1+sp}} dh'dx_Ndx'\\
&=\int_{B(0,1)} \int_0^{\eta(x')/8} \int_{B(0,\eta(x')/2)\setminus B(0,x_N)}\\
&\tab\times\frac{\left\vert \int_{\R^{N-1}} \phi(v')\left(g(x'+h'+x_Nv')-g(x'+x_Nv')\right)dv'\right\vert^p}{\Vert h'\Vert^{N-1+sp}} dh'dx_Ndx'\\
&\leq C\int_{\mathbb R^{N-1}} \int_0^{\eta(x')/8} \int_{B(0,\eta(x')/2)\setminus B(0,x_N)} \int_{B(0,1)}\\
&\tab\times\frac{\left\vert g(x'+h'+x_Nv')-g(x'+x_Nv')\right\vert^p}{\Vert h'\Vert^{N-1+sp}} dv'dh'dx_Ndx'\\
&\leq C\int_{\mathbb R^{N-1}} \int_0^{\infty} \int_{B(0,\eta(x')/2)\setminus B(0,x_N)} \int_{B(0,1)} \chi_{(0,\eta(x')/8)}(x_N)\\
&\tab\times\frac{\left\vert g(x'+h'+x_Nv')-g(x'+x_Nv')\right\vert^p}{\Vert h'\Vert^{N-1+sp}} dv'dh'dx_Ndx'.
\end{align*}
We apply a change of variables $X'=x'+x_Nv'$ to get:
\begin{align*} 
&\leq C\int_{\mathbb R^{N-1}} \int_0^{\infty} \int_{B(0,\eta(X'-x_Nv')/2)\setminus B(0,x_N)} \int_{B(0,1)} \chi_{(0,\eta(X'-x_Nv')/8)}(x_N)\\
&\tab\times\frac{\left\vert g(X'+h')-g(X')\right\vert^p}{\Vert h'\Vert^{N-1+sp}} dv'dh'dx_NdX'.
\end{align*}
We note that $\eta(X'-x_Nv')\leq \eta(X')+x_N$ since $\eta$ is $1$-Lipschitz. Since $x_N\leq \eta(X'-x_Nv')/8$, we get that $\eta(X'-x_Nv')\leq \frac{8}{7}\eta(X')$. We use this to get:
\begin{align*} 
&\leq C\int_{\mathbb R^{N-1}} \int_0^{\infty} \int_{B(0,4\eta(X')/7)\setminus B(0,x_N)} \int_{B(0,1)} \chi_{(0,4\eta(X')/7)}(x_N)\frac{\left\vert g(X'+h')-g(X')\right\vert^p}{\Vert h'\Vert^{N-1+sp}} dv'dh'dx_NdX'\\
&\leq C\int_{\mathbb R^{N-1}} \int_0^{4\eta(X')/7} \int_{B(0,4\eta(X')/7)\setminus B(0,x_N)} \frac{\left\vert g(X'+h')-g(X')\right\vert^p}{\Vert h'\Vert^{N-1+sp}} dh'dx_NdX'\\
&= C\int_{\mathbb R^{N-1}}\int_{B(0,4\eta(X')/7)}  \int_0^{\Vert h'\Vert} \frac{\left\vert g(X'+h')-g(X')\right\vert^p}{\Vert h'\Vert^{N-1+sp}} dh'dx_NdX'\\
&= C\int_{\mathbb R^{N-1}}\int_{B(0,4\eta(X')/7)} \frac{\left\vert g(X'+h')-g(X')\right\vert^p}{\Vert h'\Vert^{N-2+sp}} dh'dx_NdX'.
\end{align*}
Let $Y'=X'+h'$.
\begin{align*}
&= C\int_{\mathbb R^{N-1}}\int_{B(X',\eta(X')/2)} \frac{\left\vert g(Y')-g(X')\right\vert^p}{\Vert Y'-X'\Vert^{N-2+sp}} dY'dx_NdX'\\
&\tab+ C\int_{\mathbb R^{N-1}}\int_{B(X',4\eta(X')/7)\setminus B(X',\eta(X')/2)} \frac{\left\vert g(Y')-g(X')\right\vert^p}{\Vert Y'-X'\Vert^{N-2+sp}} dY'dx_NdX'.
\end{align*}
The hypothesis bounds the first term. For the second term:
\begin{align*}
&=C\int_{\mathbb R^{N-1}}\int_{B(X',4\eta(X')/7)\setminus B(X',\eta(X')/2)} \Vert Y'-X'\Vert^2\frac{\left\vert g(Y')-g(X')\right\vert^p}{\Vert Y'-X'\Vert^{N+sp}} dY'dx_NdX'.
\end{align*}
As we have: $$\Vert Y'-X'\Vert \leq 4\eta(X')/7$$
And:
$$\Vert Y'-X'\Vert \leq 4(\eta (Y')+\Vert Y'-X'\Vert)/7$$
$$\frac{3}{4}\Vert Y'-X'\Vert \leq \eta(Y').$$
Therefore, we get a bound:
\begin{align*}
&\leq C\int_{\mathbb R^{N-1}}\int_{B(X',4\eta(X')/7)\setminus B(X',\eta(X')/2)}\eta(X')\eta(Y')\frac{\left\vert g(Y')-g(X')\right\vert^p}{\Vert Y'-X'\Vert^{N+sp}} dY'dx_NdX'.
\end{align*}
Now we consider:
\begin{align*}
\mathcal A_{1,1}=\int_{\mathbb R^{N-1}} \int_0^{\eta(x')/8} \int_{B(0,x_N)} \frac{|u(x'+h',x_N)-u(x',x_N)|^p}{\Vert h'\Vert^{N-1+sp}} dh'dx_Ndx'
\end{align*}
Using lemma \ref{lem:lateralbound}, we get:
\begin{align*}
&\leq \int_{\mathbb R^{N-1}} \int_0^{\eta(x')/8} \int_{B(0,x_N)} \frac{|u(x'+h',x_N)-u(x',x_N)|^p}{\Vert h'\Vert^{N-1+sp}} dh'dx_Ndx'\\
&\leq C\int_{\mathbb R^{N-1}} \int_{B(0,\eta(x')/8)} \int_{\Vert h'\Vert}^{\frac{\eta(x')}{8}} \frac{1}{\Vert h'\Vert^{N-p+sp}x_N^{p+N-1}}\\
&\tab\times\int_0^{\Vert h'\Vert} \int_{B(x'+r\hat h',x_N)} |g(y')-g(x'+r\hat h')|^pdy'drdx_Ndh'dx'\\
&= C\int_{\mathbb R^{N-1}} \int_{\mathbb R^{N-1}} \int_{\Vert h'\Vert}^{\infty} \frac{1}{\Vert h'\Vert^{N-p+sp}x_N^{p+N-1}}\int_0^{\Vert h'\Vert} \int_{B(x'+r\hat h',x_N)} \chi_{(0,\eta(x')/8)}(x_N)\\
&\tab\times\chi_{B(0,\eta(x')/8)}(h')|g(y')-g(x'+r\hat h')|^pdy'drdx_Ndh'dx'.
\end{align*}
Let us use the change of variables $z'=x'+r\hat h'$. We then get:
\begin{align*}
&= C\int_{\mathbb R^{N-1}} \int_{\mathbb R^{N-1}} \int_{\Vert h'\Vert}^{\infty} \frac{1}{\Vert h'\Vert^{N-p+sp}x_N^{p+N-1}}\int_0^{\Vert h'\Vert} \int_{B(z',x_N)} \chi_{(0,\eta(z'-r\hat h')/8)}(x_N)\\
&\tab\times\chi_{B(0,\eta(z'-r\hat h')/8)}(h')|g(y')-g(z')|^pdy'drdx_Ndh'dz'.
\end{align*}
Using the fact that $\eta$ is $1$-Lipschitz and the indicator function on $h'$, we have:
$$\eta(z'-rh')\leq \eta(z')+r\leq \eta(z')+\Vert h'\Vert\leq \eta(z')+\eta(z'-rh')/8.$$
Therefore, we get:
\begin{align*}
&\leq C\int_{\mathbb R^{N-1}} \int_{\mathbb R^{N-1}} \int_{\Vert h'\Vert}^{\infty} \frac{1}{\Vert h'\Vert^{N-p+sp}x_N^{p+N-1}}\int_0^{\Vert h'\Vert} \int_{B(z',x_N)} \chi_{(0,\eta(z')/7)}(x_N)\\
&\tab\times\chi_{B(0,\eta(z')/7)}(h')|g(y')-g(z')|^pdy'drdx_Ndh'dz'\\
&= C\int_{\mathbb R^{N-1}} \int_{B(0,\eta(z')/7)} \int_{\Vert h'\Vert}^{\eta(z')/7} \frac{1}{\Vert h'\Vert^{N-p+sp}x_N^{p+N-1}}\int_0^{\Vert h'\Vert} \int_{B(z',x_N)} |g(y')-g(z')|^pdy'drdx_Ndh'dz'\\
&= C\int_{\mathbb R^{N-1}} \int_{B(z',\eta(z')/7)} |g(y')-g(z')|^p \int_{\Vert z'-y'\Vert}^{\eta(z')/7}\int_{B(0,x_N)} \frac{1}{\Vert h'\Vert^{N-p-1+sp}x_N^{p+N-1}}dh'dx_Ndy'dz' \\
&= C\int_{\mathbb R^{N-1}} \int_{B(z',\eta(z')/7)} \frac{|g(y')-g(z')|^p}{\Vert y'-z'\Vert^{N-2+sp}}dy'dz'.
\end{align*}

We now look at $\mathcal A_2$:
\begin{align*}
\mathcal A_2&=\int_{\mathbb R^{N-1}} \int_0^{\eta(x')/4} \int_0^{\eta(x')/4} \frac{|u(x',x_N)-u(x',y_N)|^p}{|x_N-y_N|^{1+sp}}dy_Ndx_Ndx'\\
&=\int_{\mathbb R^{N-1}} \int_0^{\eta(x')/4} \int_0^{\eta(x')/4} \frac{|u(x',x_N+t)-u(x',x_N)|^p}{t^{1+sp}}dtdx_Ndx'.
\end{align*}
By lemma \ref{verticalbound}, we have:
$$\leq C\int_{\mathbb R^{N-1}} \int_0^{\eta(x')/2} x_N^{p-sp} |\partial_N u(x',x_N)|^pdx_Ndx'.$$
Now, we use the inequality:
$$\vert \partial_N u(x)\vert \leq C \frac{1}{x_N^N} \int_{B(x',x_N)}|g(y')-g(x')|dy'$$
To get:
$$\leq C \int_{\mathbb R^{N-1}} \int_0^{\eta(x')/2}  x_N^{p-sp}\left(\int_{B(x',x_N)} \frac{1}{x_N^{N}} |g(y')-g(x')|dy'\right)^pdx_Ndx.'$$
By H\"older's inequality:
\begin{align*}
&\leq C \int_{\mathbb R^{N-1}} \int_0^{\eta(x')/2}  x_N^{p-sp-Np+(N-1)(p-1)}\int_{B(x',x_N)} |g(y')-g(x')|^pdy'dx_Ndx'\\
&=C\int_{\mathbb R^{N-1}} \int_{B(x',\eta(x')/2)} |g(y')-g(x')|^p \int_{\Vert x'-y'\Vert}^{\eta(x')/2} \frac{1}{x_N^{N-1+sp}}dx_Ndy'dx'\\
&\leq C\int_{\mathbb R^{N-1}}\int_{B(x',\eta(x')/2)} \frac{|g(x')-g(y')|^p}{\Vert x'-y'\Vert^{sp}}dy'dx'.
\end{align*}
\end{proof}
Now, we wish to prove that we can create a function given that the other side is $0$. Namely:
\begin{lemma}\label{generalonesidezero}
Let $g\in L^p_{\text{loc}}(\mathbb R^{N-1})$ such that:
\begin{align*}
\int_{\mathbb R^{N-1}}\int_{B(x',\eta(x')/2)} \frac{|g(x')-g(y')|^p}{\Vert x'-y'\Vert^{N-2+sp}} dy'dx'&<\infty\\
\int_{\mathbb R^{N-1}}\int_{\mathbb R^{N-1}\setminus B(x',\eta(x')/2)} \eta(x')\eta(y')\frac{|g(x')-g(y')|^p}{\Vert x'-y'\Vert^{N+sp}} dy'dx'&<\infty\\
\int_{\mathbb R^{N-1}} \frac{|g(x')|^p}{\eta(x')^{sp-1}} dx&<\infty
\end{align*}
Then there is $u\in \dot W^{s,p}(\Omega)$ such that:
\begin{align*}
\Tr(u)(x',0)&=g(x')\\
\Tr(u)(x',\eta(x'))&=0
\end{align*}
Furthermore, there is $C>0$ such that
$$C|u|_{\dot W^{s,p}(\Omega)}\leq |g|_{\tilde{W}^{s-1/p,p}_{\leq \eta}(\R^{N-1})}+|g|_{\tilde{W}^{s+1/p,p}_{\geq \eta}(\R^{N-1})}+\leq \Vert g/\eta^{s-1/p}\Vert_{L^p(\R^{N-1})}$$
\end{lemma}
\begin{proof}
We find $u_0$ such that $u_0\in \dot W^{s,p}(\Omega)$ and $\Tr(u_0)(x',0)=g(x')$ using lemma~\ref{lem:onesidegeneral}. Let $\phi$ be a $C^\infty$ function such that, for $x\leq 0$, $\phi(x)=1$ and, for $x\geq 1$, $\phi(x)=0$. We then let:
$$u(x',x_N)=\phi(x_N/\eta(x')) u_0(x',x_N)$$
We would like to show that this is still in $\dot W^{s,p}(\Omega)$. We have:
\begin{align*}
\int_{\Omega} \int_{\Omega} \frac{|u(x)-u(y)|^p}{\Vert x-y\Vert^{N+sp}}dydx&=\int_{\Omega}\int_{B(x',\eta(x')/2)} \int_0^{\eta(y')} \frac{|u(x)-u(y)|^p}{\Vert x-y\Vert^{N+sp}} dy_Ndy'dx\\
&\tab+\int_{\Omega}\int_{\mathbb R^{N-1}\setminus B(x',\eta(x')/2)} \int_0^{\eta(y')} \frac{|u(x)-u(y)|^p}{\Vert x-y\Vert^{N+sp}} dy_Ndy'dx\\
&=\mathcal A+\mathcal B.
\end{align*}
We first consider:
\begin{align*}
\mathcal B&=\int_{\mathbb R^{N-1}}\int_0^{\eta(x')}\int_{\mathbb R^{N-1}\setminus B(x',\eta(x')/2)} \int_0^{\eta(y')} \frac{|u(x)-u(y)|^p}{\Vert x'-y'\Vert^{N+sp}} dy_Ndy'dx_Ndx'\\
&\leq C\int_{\mathbb R^{N-1}}\int_0^{\eta(x')}\int_{\mathbb R^{N-1}\setminus B(x',\eta(x')/2)} \int_0^{\eta(y')} \frac{|u(x)|^p}{\Vert x'-y'\Vert^{N+sp}} dy_Ndy'dx_Ndx'\\
&\tab+C\int_{\mathbb R^{N-1}}\int_0^{\eta(x')}\int_{\mathbb R^{N-1}\setminus B(x',\eta(x')/2)} \int_0^{\eta(y')} \frac{|u(y)|^p}{\Vert x'-y'\Vert^{N+sp}} dy_Ndy'dx_Ndx'\\
&=:C\mathcal B_1+C\mathcal B_2.
\end{align*}
We consider:
\begin{align}
\mathcal B_1&=\int_{\mathbb R^{N-1}}\int_0^{\eta(x')}\int_{\mathbb R^{N-1}\setminus B(x',\eta(x')/2)} \int_0^{\eta(y')} \frac{|u(x)|^p}{\Vert x'-y'\Vert^{N+sp}} dy_Ndy'dx_Ndx'\label{B1math}\\
&\leq \int_{\mathbb R^{N-1}}\int_0^{\eta(x')}\int_{\mathbb R^{N-1}\setminus B(x',\eta(x')/2)} \int_0^{\eta(x')+\Vert x'-y'\Vert} \frac{|u(x)|^p}{\Vert x'-y'\Vert^{N+sp}} dy_Ndy'dx_Ndx'\nonumber\\
&= \int_{\mathbb R^{N-1}}\int_0^{\eta(x')}\int_{\mathbb R^{N-1}\setminus B(x',\eta(x')/2)} (\eta(x')+\Vert x'-y'\Vert) \frac{|u(x)|^p}{\Vert x'-y'\Vert^{N+sp}} dy'dx_Ndx'\nonumber\\
&\leq C\int_{\mathbb R^{N-1}}\int_0^{\eta(x')}\int_{\mathbb R^{N-1}\setminus B(x',\eta(x')/2)} \Vert x'-y'\Vert \frac{|u(x)|^p}{\Vert x'-y'\Vert^{N+sp}} dy'dx_Ndx'\nonumber\\
&\leq C\int_{\mathbb R^{N-1}}\int_0^{\eta(x')} |u(x)|^p\int_{\mathbb R^{N-1}\setminus B(x',\eta(x')/2)} \frac{1}{\Vert x'-y'\Vert^{N-1+sp}} dy'dx_Ndx'\nonumber\\
&= C\int_{\mathbb R^{N-1}}\int_0^{\eta(x')} \frac{|u(x)|^p}{\eta(x')^{sp}} dx_Ndx'.\nonumber
\end{align}
Since $\psi\leq 1$, we get
\begin{align}
&\leq C\int_{\mathbb R^{N-1}}\int_0^{\eta(x')} \frac{|u_0(x)|^p}{\eta(x')^{sp}} dx_Ndx' \label{something}\\
&\leq C\int_{\mathbb R^{N-1}}\int_0^{\eta(x')}\int_{B(0,1)} \frac{|g(x'+\frac{x_N}{8}w')|^p}{\eta(x')^{sp}} dw'dx_Ndx'\nonumber \\
&= C\int_{\mathbb R^{N-1}}\int_0^{\infty}\int_{B(0,1)} \chi_{(0,\eta(X'-\frac{x_N}{8}w'))}(x_N) \frac{|g(X')|^p}{\eta(X'-\frac{x_N}{8}w')^{sp}} dw'dx_NdX'.\nonumber 
\end{align}
Where we used the change of variables $X'=x'+\frac{x_N}{8}w'$. We note:
$$\eta(X')-\frac{\eta(X'-\frac{x_N}{8}w')}{8}\leq \eta(X')-\frac{x_N}{8}\leq \eta(X'-\frac{x_N}{8}w')\leq \eta(X')+\frac{x_N}{8}\leq \eta(X')+\frac{\eta(X'-\frac{x_N}{8}w')}{8}.$$
Therefore:
$$\frac{8}{9}\eta(X')\leq \eta(X'-\frac{x_N}{8}w')\leq \frac{8}{7}\eta(X').$$
So we get:
\begin{align*}
&\leq C\int_{\mathbb R^{N-1}}\int_0^{8\eta(X')/7}\int_{B(0,1)} \frac{|g(X')|^p}{(8\eta(X')/9)^{sp}} dw'dx_NdX'\\
&\leq C\int_{\mathbb R^{N-1}} \frac{8}{7}\eta(X') \frac{|g(X')|^p}{(8\eta(X')/9)^{sp}} dX'\\
&\leq C\int_{\mathbb R^{N-1}} \frac{|g(X')|^p}{\eta(X')^{sp-1}} dX'.
\end{align*}
Next, we consider: 
\begin{align*}
\mathcal B_2&=\int_{\mathbb R^{N-1}}\int_0^{\eta(x')}\int_{\mathbb R^{N-1}\setminus B(x',\eta(x')/2)} \int_0^{\eta(y')} \frac{|u(y)|^p}{\Vert x'-y'\Vert^{N+sp}} dy_Ndy'dx_Ndx'\\
&=\int_{\mathbb R^{N-1}}\int_0^{\eta(x')}\int_{\mathbb R^{N-1}} \chi_{(\eta(x')/2,\infty)}(\Vert x'-y'\Vert) \int_0^{\eta(y')} \frac{|u(y)|^p}{\Vert x'-y'\Vert^{N+sp}} dy_Ndy'dx_Ndx'.
\end{align*}
We note that if
$$\Vert x'-y'\Vert \geq \eta(x')/2,$$
then
$$\Vert x'-y'\Vert \geq \eta(y')/2-\Vert x'-y'\Vert/2,$$
so
$$\Vert x'-y'\Vert \geq \eta(y')/3.$$
which gets us:
\begin{align*}
\leq C\int_{\mathbb R^{N-1}}\int_0^{\eta(y')}\int_{\mathbb R^{N-1}\setminus B(y',\eta(y')/3)} \int_0^{\eta(x')} \frac{|u(y)|^p}{\Vert x'-y'\Vert^{N+sp}} dx_Ndx'dy_Ndy'.
\end{align*}
We then repeat the computations done in (\ref{B1math}). 
We then consider, by lemma \ref{lem:lateralbound}
\begin{align*}
\mathcal A&\leq C\int_{\mathbb R^{N-1}} \int_0^{\eta(x')} \int_{B(x',\eta(x')/2)} \frac{|u(x',x_N)-u(y',x_N)|^p}{\Vert x'-y'\Vert^{N-1+sp}}dy'dx_Ndx'\\
&\tab+C\int_{\mathbb R^{N-1}} \int_0^{2\eta(x')} \int_0^{2\eta(x')} \frac{|u(x',x_N)-u(x',y_N)|^p}{|x_N-y_N|^{1+sp}}dy'dx_Ndx'\\
&=:C\mathcal A_1+C\mathcal A_2.
\end{align*}
We have:
\begin{align*}
\mathcal A_1&=\int_{\mathbb R^{N-1}} \int_0^{\eta(x')} \int_{B(x',\eta(x')/2)} \frac{|\psi(x_N/\eta(x'))u_0(x',x_N)-\psi(x_N/\eta(y'))u_0(y',x_N)|^p}{\Vert x'-y'\Vert^{N-1+sp}}dy'dx_Ndx'\\
&\leq C\int_{\mathbb R^{N-1}} \int_0^{\eta(x')} \int_{B(x',\eta(x')/2)} \frac{|\psi(x_N/\eta(x'))u_0(x',x_N)-\psi(x_N/\eta(y'))u_0(x',x_N)|^p}{\Vert x'-y'\Vert^{N-1+sp}}dy'dx_Ndx'\\
&\tab+C\int_{\mathbb R^{N-1}} \int_0^{\eta(x')} \int_{B(x',\eta(x')/2)} \frac{|\psi(x_N/\eta(y'))u_0(x',x_N)-\psi(x_N/\eta(y'))u_0(y',x_N)|^p}{\Vert x'-y'\Vert^{N-1+sp}}dy'dx_Ndx'\\
&=:\mathcal A_{1,1}+\mathcal A_{1,2}.
\end{align*}
We have:
\begin{align*}
\mathcal A_{1,1}&=\int_{\mathbb R^{N-1}} \int_0^{\eta(x')} \int_{B(x',\eta(x')/2)} \frac{|u_0(x',x_N)|^p\cdot |\psi(x_N/\eta(x'))-\psi(x_N/\eta(y'))|^p}{\Vert x'-y'\Vert^{N-1+sp}}dy'dx_Ndx'\\
\end{align*}
We have $\psi$ is Lipschitz, so we get:
\begin{align*}
&\leq C\int_{\mathbb R^{N-1}} \int_0^{\eta(x')} \int_{B(x',\eta(x')/2)} \frac{|u_0(x',x_N)|^p\cdot |(x_N/\eta(x'))-(x_N/\eta(y'))|^p}{\Vert x'-y'\Vert^{N-1+sp}}dy'dx_Ndx'\\
&\leq C\int_{\mathbb R^{N-1}} \int_0^{\eta(x')} \int_{B(x',\eta(x')/2)} \frac{|u_0(x',x_N)|^p\cdot x_N^p\cdot |\eta(y')-\eta(x')|^p}{\eta(x')^p\eta(y')^p\Vert x'-y'\Vert^{N-1+sp}}dy'dx_Ndx'.
\end{align*}
Note that we have the following properties:
$$x_N\leq \eta(x'),$$
$$|\eta(y')-\eta(x')|\leq \Vert x'-y'\Vert,$$
$$\eta(y')\geq \eta(x')/2.$$
Therefore, we get:
\begin{align*}
&\leq C\int_{\mathbb R^{N-1}} \int_0^{\eta(x')} \int_{B(x',\eta(x')/2)} \frac{|u_0(x',x_N)|^p}{\eta(x')^p\Vert x'-y'\Vert^{N-1-p+sp}}dy'dx_Ndx'\\
&=C\int_{\mathbb R^{N-1}} \int_0^{\eta(x')} \frac{|u_0(x',x_N)|^p}{\eta(x')^p} \int_{B(x',\eta(x')/2)} \frac{1}{\Vert x'-y'\Vert^{N-1-p+sp}}dy'dx_Ndx'\\
&=C\int_{\mathbb R^{N-1}} \int_0^{\eta(x')} \frac{|u_0(x',x_N)|^p}{\eta(x')^{sp}}dx_Ndx'.
\end{align*}
We use (\ref{something}) to bound this. We now consider: 
\begin{align*}
\mathcal A_{1,2}&=\int_{\mathbb R^{N-1}} \int_0^{\eta(x')} \int_{B(x',\eta(x')/2)} \frac{|\psi(x_N/\eta(y'))|^p\cdot |u_0(x',x_N)-u_0(y',x_N)|^p}{\Vert x'-y'\Vert^{N-1+sp}}dy'dx_Ndx'\\
&\leq \int_{\mathbb R^{N-1}} \int_0^{\eta(x')} \int_{B(x',\eta(x')/2)} \frac{|u_0(x',x_N)-u_0(y',x_N)|^p}{\Vert x'-y'\Vert^{N-1+sp}}dy'dx_Ndx'.
\end{align*}
We already proved this is bounded in \ref{thm:equivdomaingeneral}.  Now, we consider:
\begin{align*}
\mathcal A_2&=\int_{\mathbb R^{N-1}} \int_0^{2\eta(x')} \int_0^{2\eta(x')} \frac{|u(x',x_N)-u(x',y_N)|^p}{|x_N-y_N|^{1+sp}}dy_Ndx_Ndx'\\
&\leq \int_{\mathbb R^{N-1}} \int_0^{2\eta(x')} \int_0^{2\eta(x')} \frac{|\psi(x_N/\eta(x'))u_0(x',x_N)-\psi(y_N/\eta(x'))u_0(x',y_N)|^p}{|x_N-y_N|^{1+sp}}dy_Ndx_Ndx'\\
&\leq C\int_{\mathbb R^{N-1}} \int_0^{2\eta(x')} \int_0^{2\eta(x')} \frac{|\psi(x_N/\eta(x'))u_0(x',x_N)-\psi(x_N/\eta(x'))u_0(x',y_N)|^p}{|x_N-y_N|^{1+sp}}dy_Ndx_Ndx'\\
&\tab+C\int_{\mathbb R^{N-1}} \int_0^{2\eta(x')} \int_0^{2\eta(x')} \frac{|\psi(x_N/\eta(x'))u_0(x',y_N)-\psi(y_N/\eta(x'))u_0(x',y_N)|^p}{|x_N-y_N|^{1+sp}}dy_Ndx_Ndx'\\
&=:C\mathcal A_{2,1}+C\mathcal A_{2,2}.
\end{align*}
We have:
\begin{align*}
\mathcal A_{2,1}&=\int_{\mathbb R^{N-1}} \int_0^{2\eta(x')} \int_0^{2\eta(x')} \frac{|\psi(x_N/\eta(x'))u_0(x',x_N)-\psi(x_N/\eta(x'))u_0(x',y_N)|^p}{|x_N-y_N|^{1+sp}}dy_Ndx_Ndx'\\
&=\int_{\mathbb R^{N-1}} \int_0^{2\eta(x')} \int_0^{2\eta(x')} \frac{|\psi(x_N/\eta(x'))|^p\cdot |u_0(x',x_N)-u_0(x',y_N)|^p}{|x_N-y_N|^{1+sp}}dy_Ndx_Ndx'\\
&\leq \int_{\mathbb R^{N-1}} \int_0^{2\eta(x')} \int_0^{2\eta(x')} \frac{ |u_0(x',x_N)-u_0(x',y_N)|^p}{|x_N-y_N|^{1+sp}}dy_Ndx_Ndx'.
\end{align*}
Next, we have, using the Lipschitz continuity of $\psi$:
\begin{align*}
\mathcal A_{2,2}&=\int_{\mathbb R^{N-1}} \int_0^{2\eta(x')} \int_0^{2\eta(x')} \frac{|\psi(x_N/\eta(x'))u_0(x',y_N)-\psi(y_N/\eta(x'))u_0(x',y_N)|^p}{|x_N-y_N|^{1+sp}}dy_Ndx_Ndx'\\
&=\int_{\mathbb R^{N-1}} \int_0^{2\eta(x')} \int_0^{2\eta(x')} \frac{|u_0(x',y_N)|^p\cdot |\psi(x_N/\eta(x'))-\psi(y_N/\eta(x'))|^p}{|x_N-y_N|^{1+sp}}dy_Ndx_Ndx'\\
&\leq\int_{\mathbb R^{N-1}} \int_0^{2\eta(x')} \int_0^{2\eta(x')} \frac{|u_0(x',x_N)|^p\cdot |x_N-y_N|^p}{\eta(x')^p|x_N-y_N|^{1+sp}}dy_Ndx_Ndx'\\
&=\int_{\mathbb R^{N-1}} \int_0^{2\eta(x')}\frac{|u_0(x',y_N)|^p}{\eta(x')^p} \int_0^{2\eta(x')} \frac{1}{|x_N-y_N|^{1-p+sp}}dx_Ndy_Ndx'\\
&\leq\int_{\mathbb R^{N-1}} \int_0^{2\eta(x')}\frac{|u_0(x',y_N)|^p}{\eta(x')^p} \int_0^{2\eta(x')} \frac{1}{|x_N-y_N|^{1-p+sp}}dx_Ndy_Ndx'\\
&\leq\int_{\mathbb R^{N-1}} \int_0^{2\eta(x')}\frac{|u_0(x',y_N)|^p}{\eta(x')^{sp}} dy_Ndx'.
\end{align*}
Apply (\ref{something}) to bound. 
\end{proof}

%% file: sections/5.Fourier_Transforms.tex
\section{Comparing Spaces}
\subfile{5.Fourier_Transforms/5.Fourier_Transforms}

%% file: sections/5.Fourier_Transforms/5.Fourier_Transforms.tex
While we introduced the space $W^{s,p}_{\geq \eta}(\Omega)$ to define our trace, it is still unclear how this space relates to previously researched screened Sobolev spaces. 
\par In this section, we discuss the relations of these spaces for $\eta=1$ on $\R$. We have the containment
$$\dot W^{s-1/p,p}(\R)\subseteq  \tilde{W}^{s-1/p,p}_{\leq 1}(\R)\cap \tilde{W}^{s+1/p}_{\geq 1}(\R)\subseteq \tilde{W}^{s-1/p}_{\leq 1}(\R).$$
We show for certain values of $s,p$ that the containment are strict using examples of functions. Furthermore, we will show, using the Fourier equivalence of these seminorms, the equality of these spaces for some values of $s,p$. 
\subsection{Strict Containment of the Spaces through Examples}
We hope to show that, for specific values of $s,p$, that
$$\dot W^{s-1/p,p}(\R)\subsetneqq \tilde{W}^{s-1/p,p}_{\leq 1}(\R)\cap \tilde{W}^{s+1/p,p}_{\geq 1}(\R)\subsetneqq \tilde{W}^{s-1/p,p}_{\leq 1}(\R).$$
Thus, we will show that $\tilde{W}^{s+1/p,p}_{\geq 1}(\Omega)$, in some cases, is necessary to describe the trace space. 
\begin{example}
Let $1\leq p<\infty$ and $0<s<1$ such that $1<sp<2$. Then
$$\chi_{[0,\infty)}\in \left(W^{s-1/p,p}_{\leq 1}(\R)\cap W^{s+1/p,p}_{\geq 1}(\R)\right)\setminus W^{s-1/p,p}(\R)$$
\end{example}
\begin{proof}
We consider
\begin{align*}
|\chi_{[0,\infty)}|_{\dot W^{s-1/p,p}(\R)}^p&=\int_\R\int_\R \frac{\left\vert \chi_{[0,\infty)}(x)-\chi_{[0,\infty)}(y)\right\vert^p}{|x-y|^{sp}}dydx\\
&=2\int_0^\infty \int_{-\infty}^0 \frac{1}{|x-y|^{sp}}dydx\\
&=2\int_0^\infty \int_x^\infty \frac{1}{r^{sp}}drdx\\
&=\frac{2}{sp-1}\int_0^\infty \frac{1}{x^{sp-1}}dx.
\end{align*}
Note that this is unbounded. Therefore, instead we consider
\begin{align*}
|\chi_{[0,\infty)}|_{\tilde{W}^{s-1/p,p}_{\leq 1}(\R)}^p&=\int_\R\int_{B(x,1)} \frac{\left\vert \chi_{[0,\infty)}(x)-\chi_{[0,\infty)}(y)\right\vert^p}{|x-y|^{sp}}dydx\\&=2\int_0^\infty \int_{B(x,1)} \frac{\left\vert \chi_{[0,\infty)}(x)-\chi_{[0,\infty)}(y)\right\vert^p}{|x-y|^{sp}}dydx\\
&=2\int_0^1 \int_x^1 \frac{1}{r^{sp}}drdx\\
&=\frac{2}{sp-1}\int_0^1 \frac{1}{x^{sp-1}}-1dx.
\end{align*}
Note that if $sp<2$, then this is bounded. Meanwhile,
\begin{align*}
|\chi_{[0,\infty)}|_{\tilde{W}^{s+1/p,p}_{\geq 1}(\R)}^p&=\int_\R \int_{\R\setminus B(x,1)} \frac{\left\vert \chi_{[0,\infty)}(x)-\chi_{[0,\infty)}(y)\right\vert^p}{|x-y|^{2+sp}}dydx \\
&= 2\int_0^\infty \int_{\R\setminus B(x,1)} \frac{\left\vert \chi_{[0,\infty)}(x)-\chi_{[0,\infty)}(y)\right\vert^p}{|x-y|^{2+sp}}dydx\\
&= 2\int_0^\infty \int_{\max(1,x)}^\infty \frac{1}{r^{2+sp}}drdx\\
&= 2\int_0^1 \int_{1}^\infty \frac{1}{r^{2+sp}}drdx+2\int_1^\infty \int_{x}^\infty \frac{1}{r^{2+sp}}drdx\\
&= \frac{2}{1+sp}\int_0^1 1 dx+\frac{2}{1+sp}\int_1^\infty \frac{1}{x^{1+sp}}dx.
\end{align*}
We note that this is bounded. Therefore, with this example, we have
$$W^{s-1/p,p}(\R)\subsetneqq W^{s-1/p,p}_{\leq 1}(\R)\cap W^{s+1/p,p}_{\geq 1}(\R).$$
\end{proof}
\begin{example}
Let $0<s<1$, $2< p<\infty$ and $0<\lambda<1$ such that $1/p<s<\lambda < 1-1/p$. Furthermore, define
$$f(x)=\max(1,x^\lambda).$$
Then
$$f\in W^{s-1/p,p}_{\leq 1}(\R)\setminus W^{s+1/p,p}_{\geq 1}(\R).$$
\end{example}
\begin{proof}
First note that for $x\geq 1$, the derivative is $\lambda x^{\lambda-1}\leq \lambda$. Furthermore, the derivative for $x\geq 1$ is decreasing. Therefore, by the mean value theorem,
\begin{align*}
|f|_{W^{s-1/p,p}_{\leq 1}(\R)}^p&=\int_\R\int_{B(x,1)} \frac{\left\vert f(x)-f(y)\right\vert^p}{|x-y|^{sp}}dydx\\&=2\int_\R \int_0^1 \frac{|f(x+r)-f(x)|^p}{r^{sp}}drdx\\
&=2\int_0^1 \int_0^1 \frac{|f(x+r)-f(x)|^p}{r^{sp}}drdx+2\int_1^\infty \int_0^1 \frac{|f(x+r)-f(x)|^p}{r^{sp}}drdx\\
&\leq 2\int_0^1 \int_0^1 \frac{(\lambda r)^p}{r^{sp}}drdx+2\int_1^\infty \int_0^1 \frac{|\lambda x^{\lambda-1}r|^p}{r^{sp}}drdx\\
&\leq C \int_0^1 r^{p(1-s)}dr+C \int_1^\infty x^{p(\lambda-1)} dx\\
&<\infty
\end{align*}
On the other hand,
\begin{align*}
|f|_{W^{s+1/p,p}_{\geq 1}(\R)}^p&=\int_\R\int_{\R\setminus B(x,1)} \frac{\left\vert f(x)-f(y)\right\vert^p}{|x-y|^{2+sp}}dydx\\&\geq \int_1^\infty \int_{-\infty}^0 \frac{(x^\lambda -1)^p}{|x-y|^{2+sp}}dydx\\
&= \int_1^\infty (x^\lambda -1)^p\int_x^\infty \frac{1}{r^{2+sp}}drdx\\
&= C\int_1^\infty  \frac{(x^\lambda -1)^p}{x^{1+sp}}dx\\
\end{align*}
We note that if $\lambda p-1-sp>-1$, then this is unbounded. Since $\lambda>s$, then $$\lambda p-1-sp>sp-1-sp=-1.$$
\end{proof}
\subsection{Equivalence of Certain Spaces with Fourier Transforms}
While we have strict containment of the spaces for certain values of $s,p$, we can show that the spaces are equivalent for other values. We approach this problem by using the equivalent Fourier seminorms for $p=2$. We draw inspiration from proposition $3.4$ in \cite{di_nezza_palatucci_valdinoci_2012}. 
\begin{theorem}
Let $1/2\leq s<1$. Then $g\in \tilde{W}^{s-1/2,2}_{\leq 1}(\R)\cap \tilde W^{s+1/2,p}_{\geq 1}(\R)$ if and only if
$$\int_{|x|\leq 1/2} |x|^2\cdot |\hat g(x)|^2dx+\int_{|x|\geq 1/2} |x|^{2s-1}\cdot |\hat g(x)|^2 dx<\infty$$ 
\end{theorem}
\begin{proof}
Note that $g\in \tilde{W}^{s-1/2,2}_{\leq 1}(\R)\cap \tilde W^{s+1/2,p}_{\geq 1}(\R)$ if and only if
$$\int_{\R}\int_{B(x,1)} \frac{|g(x)-g(y)|^2}{|x-y|^{2s}}dydx+\int_{\R}\int_{\R\setminus B(x,1)} \frac{|g(x)-g(y)|^2}{|x-y|^{2+2s}}dydx<\infty.$$
We note that we can combine the integrals in the form
\begin{align*}
\int_{\R}\int_{\R} \min(1,|x-y|^2) \frac{|g(x)-g(y)|^2}{|x-y|^{2+2s}}dydx&=\int_{\R}\int_{\R} \min(1,|h|^2) \frac{|g(x)-g(x+h)|^2}{|h|^{2+2s}}dhdx
\end{align*}
Let:
$$g_h(x)=\min(1,|h|)\frac{g(x+h)-g(x)}{|h|^{1+s}}$$
Then, we have that:
\begin{align*}
\firstline\int_{\R}\int_{\R}\min(1,|h|^2)\frac{|g(x+h)-g(x)|^2}{|h|^{2+2s}}\;dx\;dh\\
&=\int_{\R}\int_{\R}|g_h(x)|^2\;dx\;dh.
\end{align*}
By Plancherel theorem:
\begin{align*}
\int_{\R}\int_{\R}| g_h(x)|^2\;dx\;dh=\int_{\R}\int_{\R}\left\vert\mathcal F\left(\min(1,|h|)\frac{g(\cdot+h)-g(\cdot)}{|h|^{1+s}}\right)(x)\right\vert^2\;dx\;dh.
\end{align*}
We have:
\begin{align*}
\firstline \int_{\R}\int_{\R}\left\vert\mathcal F\left(\min(1,|h|)\frac{g(\cdot+h)-g(\cdot)}{|h|^{1+s}}\right)(x)\right\vert^2\;dx\;dh\\
&=\int_{\R}\int_{\R}\frac{\min(1,|h|^2)}{|h|^{2+2s}}\left\vert\mathcal F\left(g(\cdot+h)-g(\cdot)\right)(x)\right\vert^2\;dx\;dh\\
&=\int_{\R}\int_{\R}\frac{\min(1,|h|^2)}{|h|^{2+2s}}\left\vert\int_{\R} e^{-2\pi ixy} (g(y+h)-g(y))\;dy\right\vert^2\;dx\;dh\\
&=\int_{\R}\int_{\R}\frac{\min(1,|h|^2)}{|h|^{2+2s}}\left\vert\int_{\R} e^{-2\pi ixy} g(y+h)dy-\int_{\R} e^{-2\pi ixy}g(y)\;dy\right\vert^2\;dx\;dh.
\end{align*}
 Performing a change of variables:
 \begin{align*}
\firstline \int_{\R}\int_{\R}\frac{\min(1,|h|^2)}{|h|^{2+2s}}\left\vert\int_{\R} e^{-2\pi ix\cdot (y-h)} g(y)dy-\int_{\R} e^{-2\pi ixy}g(y)\;dy\right\vert^2\;dx\;dh\\
&=\int_{\R}\int_{\R}\frac{\min(1,|h|^2)}{|h|^{2+2s}}\left\vert\int_{\R} (e^{-2\pi ix\cdot (y-h)} -e^{-2\pi ixy})g(y)\;dy\right\vert^2\;dx\;dh\\
&=\int_{\R}\int_{\R}\frac{\min(1,|h|^2)}{|h|^{2+2s}} \left\vert(e^{2\pi ixh} -1)\int_{\R}e^{-2\pi ixy}g(y)\;dy\right\vert^2\;dx\;dh\\
&=\int_{\R}\int_{\R}\frac{\min(1,|h|^2)}{|h|^{2+2s}} \left\vert(e^{2\pi ixh} -1)\mathcal F(g)(x)\right\vert^2\;dx\;dh\\
&=\int_{\R}\int_{\R}\frac{\min(1,|h|^2)}{|h|^{2+2s}} |\cos(2\pi (xh)+i\sin(2\pi (xh)-1|^2\left\vert\mathcal F(g)(x)\right\vert^2\;dx\;dh\\
&=\int_{\R}\int_{\R}\frac{\min(1,|h|^2)}{|h|^{2+2s}} ((\cos(2\pi (xh))-1)^2+\sin^2(2\pi (xh)))\left\vert\mathcal F(g)(x)\right\vert^2\;dx\;dh\\
&=\int_{\R}\int_{\R}\frac{\min(1,|h|^2)}{|h|^{2+2s}} (2-2\cos(2\pi (xh)))\left\vert\mathcal F(g)(x)\right\vert^2\;dx\;dh\\
&=\int_{\R}\left\vert\mathcal F(g)(x)\right\vert^2\int_{\R}\frac{\min(1,|h|^2)(2-2\cos(2\pi (xh)))}{|h|^{2+2s}} \;dh\;dx.
 \end{align*}
 We consider just the inner integral; we apply a change of variables $H=|x|h$ to get:
  \begin{align*}
\firstline \int_{\R}\frac{1}{|x|}\frac{\min(1,(|H|/|x|)^2)(2-2\cos(2\pi H))}{(|H|/|x|)^{2+2s}} \;dh\\
&=|x|^{2s+1}\int_{\R}\frac{\min(1,(|H|/|x|)^2)(2-2\cos(2\pi H))}{|H|^{2+2s}} \;dh\\
&=|x|^{2s-1}\int_{B(0,|x|)}\frac{2-2\cos(2\pi H)}{|H|^{2s}} \;dh\\
&\tab+|x|^{2s+1}\int_{\R\setminus B(0,|x|)}\frac{2-2\cos(2\pi H)}{|H|^{2+2s}} \;dh\\
&=\mathcal A+\mathcal B.
\end{align*}
We then consider these separately. We have that
\begin{align*}
\mathcal A=|x|^{2s-1}\int_{B(0,|x|)}\frac{2-2\cos(2\pi H)}{|H|^{2s}} \;dh.
\end{align*}
If $|x|\leq 1/2$, we note that $2-2\cos(2\pi h)$ has magnitude $h^2$; therefore, we can bound $cH^2<2-2\cos (2\pi H)<CH^2$ for some $c,C>0$ for $|H|\leq 1/2$. We then get:
\begin{align*}
c|x|^{2s-1}\int_{B(0,|x|)}\frac{H^2}{|H|^{2s}} \;dh\leq \mathcal A\leq C|x|^{2s-1}\int_{B(0,|x|)}\frac{H^2}{|H|^{2s}} \;dh.
\end{align*}
We can evaluate the integrals to get:
\begin{align*}
c|x|^{2} \leq \mathcal A\leq C|x|^{2s}.
\end{align*}
For $|x|\geq 1/2$, we spit the integral into:
\begin{align*}
\firstline  |x|^{2s-1}\int_{B(0,1/2)}\frac{2-2\cos(2\pi H)}{|H|^{2s}} \;dh+|x|^{2s-1}\int_{B(0,| x'|)\setminus B(0,1/2)}\frac{2-2\cos(2\pi H)}{|H|^{2s}} \;dh  \\
 &\approx C |x|^{2s-1} + |x|^{2s-1} \int_{B(0,\Vert x'\Vert)\setminus B(0,1/2)} \frac{2}{|H|^{2s}}\\
 &=C|x|^{2s-1}+C.
\end{align*}
Lastly, we consider $\mathcal B$. We have: 
\begin{align*}
&|x|^{2s+1}\int_{\R\setminus B(0,|x|)}\frac{2-2\cos(2\pi H)}{|H|^{2+2s}} \;dh.
\end{align*}
If $|x|\geq 1/2$, then:
\begin{align*}
&|x|^{2s+1}\int_{\R\setminus B(0,|x|)}\frac{2-2\cos(2\pi H)}{|H|^{2+2s}} \;dh\\
&\approx |x|^{2s+1}\int_{\R\setminus B(0,|x|)}\frac{2}{|H|^{2+2s}} \;dh\\
&= C.
\end{align*}
If $|x| \leq 1/2$, then:
\begin{align*}
&|x|^{2s+1}\int_{\R\setminus B(0,1/2)}\frac{2-2\cos(2\pi H)}{|H|^{2+2s}} \;dh+|x|^{2s+1}\int_{B(0,1/2)\setminus B(0,|x|)}\frac{2-2\cos(2\pi H)}{|H|^{2+2s}} \;dh\\
&\approx|x|^{2s+1}\int_{\R\setminus B(0,1/2)}\frac{2}{|H|^{2+2s}} \;dh+|x|^{2s+1}\int_{B(0,1/2)\setminus B(0,|x|)}\frac{|H|^2}{|H|^{2+2s}} \;dh\\
&=C|x|^{2s+1}+C|x|^2.
\end{align*}
\par Now note that if we put everything together, if $|x|\leq 1/2$, we get a term of the form:
$$C(|x|^2+|x|^{2s+1})$$
Note that since $2s+1>2$, the term is dominated by $|x|^2$. Meanwhile, for $|x|\geq 1/2$, we get:
$$C(|x|^{2s-1}+1).$$
Note that, since $2s>1$, the term is dominated by $|x|^{2s-1}$. Therefore, we get, by separating the value of $x$:
$$\int_{|x|\leq 1/2} |x|^2\cdot |\hat g(x)|^2dx+\int_{|x|\geq 1/2} |x|^{2s-1}\cdot |\hat g(x)|^2 dx.$$ 
\end{proof}
Now, we perform a similar computation for just $W^{s-1/p,p}_{\leq 1}(\R)$. 
\begin{theorem}
Let $1/2\leq s<1$. Then $g\in \tilde{W}^{s-1/2,2}_{\leq 1}(\R)$ if and only if
$$\int_{|x|\leq 1/2} |x|^2\cdot |\hat g(x)|^2dx+\int_{|x|\geq 1/2} |x|^{2s-1}\cdot |\hat g(x)|^2 dx<\infty$$ 
\end{theorem}
\begin{proof}
Note that $g\in \tilde{W}^{s-1/2,2}_{\leq 1}(\R)$ if and only if
$$\int_{\R}\int_{B(x,1)} \frac{|g(x)-g(y)|^2}{|x-y|^{2s}}dydx=\int_{\R}\int_{\R} \chi_{[0,1]}(|h|)\frac{|g(x+h)-g(x)|^2}{|h|^{2s}}dhdx<\infty.$$

Let:
$$g_h(x)=\chi_{[0,1]}(|h|)\frac{g(x+h)-g(x)}{|h|^{s}}$$
Then, we have that:
\begin{align*}
\firstline\int_{\R}\int_{\R}\chi_{[0,1]}(|h|)\frac{|g(x+h)-g(x)|^2}{|h|^{2s}}\;dx\;dh\\
&=\int_{\R}\int_{\R}|g_h(x)|^2\;dx\;dh.
\end{align*}
By Plancherel theorem:
\begin{align*}
\int_{\R}\int_{\R}| g_h(x)|^2\;dx\;dh&=\int_{\R}\int_{\R}\left\vert\mathcal F\left(\chi_{[0,1]}(|h|)\frac{g(\cdot+h)-g(\cdot)}{|h|^{s}}\right)(x)\right\vert^2\;dx\;dh\\
&=\int_{\R}\left\vert\mathcal F(g)(x)\right\vert^2\int_{\R}\frac{\chi_{[0,1]}(|h|)(2-2\cos(2\pi (xh)))}{|h|^{2s}} \;dh\;dx.
 \end{align*}
 We consider just the inner integral; we apply a change of variables $H=|x|h$ to get:
  \begin{align*}
\firstline \int_{\R}\frac{1}{|x|}\frac{\chi_{[0,1]}(|H|/|x|)(2-2\cos(2\pi H))}{(|H|/|x|)^{2s}} \;dh\\
&=|x|^{2s-1}\int_{\R}\frac{\chi_{[0,1]}(|H|/|x|)(2-2\cos(2\pi H))}{|H|^{2s}} \;dh\\
&=|x|^{2s-1}\int_{B(0,|x|)}\frac{2-2\cos(2\pi H)}{|H|^{2s}} \;dh=:\mathcal A.
\end{align*}
If $|x|\leq 1/2$, we note that $2-2\cos(2\pi h)$ has magnitude $h^2$; therefore, we can bound $cH^2<2-2\cos (2\pi H)<CH^2$ for some $c,C>0$ for $|H|\leq 1/2$. We then get:
\begin{align*}
c|x|^{2s-1}\int_{B(0,|x|)}\frac{H^2}{|H|^{2s}} \;dh\leq \mathcal A\leq C|x|^{2s-1}\int_{B(0,|x|)}\frac{H^2}{|H|^{2s}} \;dh.
\end{align*}
We can evaluate the integrals to get:
\begin{align*}
c|x|^{2} \leq \mathcal A\leq C|x|^{2s}. 
\end{align*}
For $|x|\geq 1/2$, we spit the integral into:
\begin{align*}
\firstline  |x|^{2s-1}\int_{B(0,1/2)}\frac{2-2\cos(2\pi H)}{|H|^{2s}} \;dh+|x|^{2s-1}\int_{B(0,| x'|)\setminus B(0,1/2)}\frac{2-2\cos(2\pi H)}{|H|^{2s}} \;dh  \\
 &\approx C |x|^{2s-1} + |x|^{2s-1} \int_{B(0,\Vert x'\Vert)\setminus B(0,1/2)} \frac{2}{|H|^{2s}}\\
 &=C|x|^{2s-1}+C.
\end{align*}
\par Now note that if we put everything together, if $|x|\leq 1/2$, we get a term of the form:
$$C(|x|^2).$$
Meanwhile, for $|x|\geq 1/2$, we get:
$$C(|x|^{2s-1}+1).$$
Note that, since $2s>1$, the term is dominated by $|x|^{2s-1}$. Therefore, we get, by separating the value of $x$:
$$\int_{|x|\leq 1/2} |x|^2\cdot |\hat g(x)|^2dx+\int_{|x|\geq 1/2} |x|^{2s-1}\cdot |\hat g(x)|^2 dx.$$ 
\end{proof}
Therefore, $\tilde{W}^{s-1/2,2}_{\leq 1}(\R)$ and $\tilde{W}^{s-1/2,2}_{\leq 1}(\R)\cap \tilde{W}^{s-1/2,2}_{\geq 1}(\R)$ are the same space for $1/2<s<1$ and $p=2$.